\documentclass[12pt]{article}
\usepackage{amsfonts}
\usepackage{full page,
arcs, amssymb, amscd, amsmath, graphicx, 
}
\newtheorem{thm}{Theorem}[section]
\newtheorem{defn}[thm]{Definition}
\newtheorem{prop}[thm]{Proposition}
\newtheorem{cor}[thm]{Corollary}
\newtheorem{lemma}[thm]{Lemma}
\newtheorem{rema}[thm]{Remark}

\newtheorem{assum}[thm]{Assumption}
\newtheorem{exam}[thm]{Example}
\newtheorem{data}[thm]{Data}

\newcommand{\halmos}{\rule{1ex}{1.4ex}}

\newcommand{\nn}{\nonumber \\}

 \newcommand{\res}{\mbox{\rm Res}}

\renewcommand{\hom}{\mbox{\rm Hom}}

 \newcommand{\pf}{{\it Proof.}\hspace{2ex}}
 
 \newcommand{\epfv}{\hspace*{\fill}\mbox{$\halmos$}\vspace{1em}}

\newcommand{\wt}{\mbox{\rm wt}\,}

\newcommand{\mbar}{\Big\vert}

\newcommand{\A}{\mathcal{A}}

\newcommand{\C}{\mathbb{C}}
\newcommand{\Z}{\mathbb{Z}}
\newcommand{\R}{\mathbb{R}}

\newcommand{\N}{\mathbb{N}}

\newcommand{\one}{\mathbf{1}}

\title{ {\bf Generators, spanning sets and existence of 
twisted modules for a grading-restricted vertex (super)algebra} }
\date{}
\author{Yi-Zhi Huang}

\begin{document}

\bibliographystyle{alpha}
\maketitle
\begin{abstract}
For a grading-restricted vertex superalgebra $V$ and 
an automorphism $g$ of $V$, we give a linearly independent set of generators 
of the universal lower-bounded generalized $g$-twisted $V$-module 
$\widehat{M}^{[g]}_{B}$ constructed by the author in \cite{H-const-twisted-mod}. 
We prove that there exist irreducible 
lower-bounded generalized $g$-twisted $V$-modules by showing that there exists a maximal 
proper submodule of $\widehat{M}^{[g]}_{B}$ for a one-dimensional space $M$. 
We then give several spanning sets of $\widehat{M}^{[g]}_{B}$
and discuss the relations among elements of the spanning sets. 
Assuming that $V$ is a M\"{o}bius vertex superalgebra (to make sure that 
lowest weights make sense) and that $P(V)$ (the set 
of all numbers of the form $\Re(\alpha)\in [0, 1)$ for $\alpha\in \C$ such that $e^{2\pi i \alpha}$ 
is an eigenvalue of $g$) has no accumulation point in $\R$
(to make sure that irreducible  lower-bounded generalized
$g$-twisted $V$-modules have lowest weights). 
Under suitable additional conditions, which hold when
 the twisted zero-mode algebra or the twisted Zhu's algebra 
is finite dimensional, we prove 
that there exists an irreducible grading-restricted generalized 
$g$-twisted $V$-module, which is in fact an irreducible  ordinary $g$-twisted $V$-module
when $g$ is of finite order. We also prove that 
every lower-bounded generalized module with an action of $g$ for the fixed-point subalgebra $V^{g}$
of $V$ under $g$ can be extended to a lower-bounded generalized $g$-twisted $V$-module. 
\end{abstract}

\renewcommand{\theequation}{\thesection.\arabic{equation}}
\renewcommand{\thethm}{\thesection.\arabic{thm}}
\setcounter{equation}{0}
\setcounter{thm}{0}
\section{Introduction}

In the representation theory of vertex operator algebras 
and orbifold conformal field theory, the existence of twisted modules associated to 
an automorphism of a vertex operator 
algebra has been an explicitly stated conjecture since mid 1990's. 
While there exists at least one module for a vertex operator algebra (the vertex 
operator algebra itself),  it is not obvious at all why there must be a twisted 
module for a general vertex operator algebra. 

Assuming that the vertex operator algebra is simple and $C_{2}$-cofinite and the automorphism
is of finite order, 
Dong, Li and Mason proved the existence of an
irreducible twisted module \cite{DLM2}. But no progress has been made in the 
general case for more than twenty years. 
Recently, mathematicians and 
physicists have discovered that some classes of vertex operator algebras
that are not $C_{2}$-cofinite have a very rich and 
exciting representation theory. For example,  vertex operator algebras
associated to affine Lie algebras at admissible levels are not $C_{2}$-cofinite. 
But the category of ordinary modules for such a vertex operator algebra 
still has finitely many irreducible modules and every such module is completely reducible
(see \cite{AM1} for the conjecture and a proof in the case of $\frak{sl}_{2}$ and 
\cite{A} for a proof in the general case). 
Moreover, this category has a ribbon category structure and 
in many cases even has a modular tensor category 
structure (see \cite{CHY} for a proof in the case of $\frak{sl}_{2}$ and the conjectures in the 
general case and see \cite{C}
for a proof of the rigidity in the simply-laced case).  The simple $W$-algebas
at nondegenerate admissible level can also be realized
as cosets by such a vertex operator algebra (see \cite{ACL}) and this realization 
plays an important role in the proof 
of the rigidity in the simply-laced case in \cite{C}. The characters of modules for such a vertex operator algebra
also have a certain modular invariance property (see \cite{AK}). 
It is thus important to study 
vertex operator algebras that are not $C_{2}$-cofinite. 

Twisted modules are important for the study of orbifold conformal field theories. 
See for example \cite{H-problems} for two orbifold theory conjectures on the 
associativity, commutativity and modular invariance of twisted intertwining 
operators among twisted modules  introduced in \cite{H-log-twisted-mod}
and on the $G$-crossed tensor category 
structure (see \cite{T}) on the category of twisted modules. 
On the other hand, the triplet vertex operator algebras $\mathcal{W}(p, q)$  are  
kernels of so-called screening operators (see \cite{FGST}, \cite{AM} and \cite{H-log-twisted-mod}). 
These kernels can be reinterpreted as the
fixed-point subalgebras of some vertex operator (super)algebras under the automorphisms 
obtained by exponentiating the screening operators. Note that 
the automorphism given by the exponential of the screening operator for a
triplet vertex operator algebra  is  of infinite order and not semisimple. If  every module for the fixed 
point subalgebra is a submodule of a twisted module for the larger verex operator 
(super)algebra,
we can study the representation theory of these vertex operator algebras using 
the twisted representation theory of the larger vertex operator (super)algebras. 
This connection is another sign that the study of twisted modules associated to 
nonsemisimple automorphisms of infinite orders for 
vertex operator (super)algebras that are not $C_{2}$-cofinite might lead us to 
truly deep mathematics. 

The method used in \cite{DLM2} cannot be adapted to prove the existence of twisted modules 
in the general 
case since some results on genus-one $1$-point correlation functions are used there. 
Note that $C_{2}$-cofinitess is needed in the construction
of these genus-one $1$-point correlation functions. Without the $C_{2}$-cofiniteness,
we will not be able to use such correlation functions. 
This might be one of the main reasons why there has not  been much progress on the existence conjecture
in the general case for so many years.

In the present paper, we prove several conjectures on the existence of various types of $g$-twisted $V$-modules
for a grading-restricted vertex superalgebra or a M\"{o}bius vertex algebra $V$ and an automorphism $g$ of $V$. 
The automorphism $g$ does not have to be 
of finite order. Our approach is based on
the universal lower-bounded generalized $g$-twisted $V$-modules constructed by the author
in Section 5 of \cite{H-const-twisted-mod}. In fact, although these modules are constructed explicitly, the author 
did not establish in \cite{H-const-twisted-mod} that they are not $0$.
In this paper, we give explicitly a linearly independent set of
generators for such a module. In particular, these modules are indeed not $0$ and thus 
the conjecture that there exist nonzero lower-bounded generalized $g$-twisted $V$-modules is valid. 
One immediate consequence 
of this result is a proof of another conjecture stating
that the twisted Zhu's algebra (introduced in \cite{DLM1} when
$g$ is of finite order and in \cite{HY} when $g$ is general),
or equivalently, the twisted zero-mode algebra (introduced in \cite{HY}) is not $0$.
The conjecture that there are irreducible 
lower-bounded generalized $g$-twisted $V$-modules follows easily since a standard augument 
shows that there exist maximal proper submodules of 
these nonzero universal lower-bounded generalized $g$-twisted $V$-modules.

Though we obtain a set of generators of the universal lower-bounded 
generalized $g$-twisted $V$-module $\widehat{M}^{[g]}_{B}$ constructed
in \cite{H-const-twisted-mod}, the study of spanning sets for $\widehat{M}^{[g]}_{B}$
are much more difficult. We initiate this study in this paper. 
We first prove that 
the weak commutativity for twisted fields in the assumption of 
Theorem 4.3 in \cite{H-const-twisted-mod} is in fact a consequence of the other assumptions.
Then we give several spanning sets of $\widehat{M}^{[g]}_{B}$ and discuss their relations. 

To discuss the existence of irreducible grading-restricted generalized or ordinary 
$g$-twisted $V$-modules, we need to assume that the grading-restricted vertex (super)algebra $V$
be a M\"{o}bius 
vertex (super)algebra or quasi-vertex (super)algebra (see \cite{FHL}, \cite{HLZ} and 
\cite{H-twist-vo}). Note that for a M\"{o}bius 
vertex (super)algebra $V$, a lower-bounded generalized $g$-twisted $V$-module must have a 
compatible $\mathfrak{sl}_{2}$-module structure. 
Using one of the spanning set mentioned above, we are able to construct
nonzero universal lower-bounded generalized $g$-twisted $V$-modules for a M\"{o}bius 
vertex algebra $V$ by constructing a compatible $\mathfrak{sl}_{2}$-module structure on $\widehat{M}^{[g]}_{B}$.
Let $P(V)$ be the set
of all numbers $\alpha\in [0, 1)+i \R$ for $\alpha\in \C$ such that $e^{2\pi i \alpha}$ 
is an eigenvalue of $g$.
We assume in addition that $P(V)$ has no accumulation point in $\R$ to make sure that irreducible 
 lower-bounded generalized $g$-twisted $V$-modules have lowest weights. 
Then under the conditions that the set of the real parts of the lowest weights of the irreducible 
 lower-bounded generalized $g$-twisted $V$-modules has a maximum and that 
the lowest weight subspace of an irreducible 
 lower-bounded generalized $g$-twisted $V$-module with this maximum as its weight is finite dimensional, 
we prove the conjecture that there exists an irreducible grading-restricted generalized $g$-twisted $V$-module.
In particular, when the twisted zero-mode algebra or the twisted Zhu's algebra 
is finite dimensional, we prove as a consequence the conjecture
that there exists an irreducible grading-restricted generalized 
$g$-twisted $V$-module. 
In the case that the automorphism is of finite order, the irreducible 
grading-restricted generalized $g$-twisted $V$-module is in fact an ordinary $g$-twisted $V$-module. 
Note that our existence result removes the conditions that $A_{g}$ is not $0$ in 
Theorem 9.1 in  \cite{DLM1}, the simplicity of $V$ and
$C_{2}$-cofiniteness in Theorem 9.1 in \cite{DLM2} and 
weakens the conditions that $V$ has a conformal vector and 
$g$ is of finite order in these results in \cite{DLM1} and \cite{DLM2}.

Under some strong conditions on $V$
and the fixed point subalgebra $V^{G}$  under a finite group $G$ of automorphisms of $V$, including in particular, 
the conditions that both $V$ and $V^{G}$ are $C_{2}$-cofinite and reductive in the sense
that every $\N$-gradable modules is a direct sum of irreducible module, Dong, Ren and Xu proved in \cite{DRX}
that every irreducible $V^{G}$-module can be extended to an irreducible $g$-twisted $V$-module for 
some $g\in G$. 
Using the universal lower-bounded generalized twisted modules
constructed in Section 5 of \cite{H-const-twisted-mod},
we also prove another existence result 
that every lower-bounded generalized module with an action of $g$ 
for the fixed-point subalgebra $V^{g}$
of $V$ under $g$ can be extended to a lower-bounded generalized $g$-twisted $V$-module. 
This result allows us to construct and study lower-bounded generalized $V^{g}$-modules using 
lower-bounded generalized $g$-twisted $V$-modules.
Note that we do not assume any conditions on $V$ or $V^{g}$ and
the condition in our result is in fact necessary.

In this paper, we fix a grading-restricted vertex superalgebra  $V$
and an automorphism $g$ of $V$. In Section 6, $V$ 
is a fixed M\"{o}bius vertex superalgebra.  For the precise meaning of 
these notions of vertex algebra and various notions of modules used in this 
paper, see \cite{H-twist-vo}. We shall call a grading-restricted generalized 
$g$-twisted $V$-module on which $L(0)$ acts semisimply 
an ordinary $g$-twisted $V$-module. 
 Note that as in \cite{H-2-const},
\cite{H-twist-vo} and \cite{H-const-twisted-mod}, grading-restricted vertex algebras
and M\"{o}bius vertex algebras are special cases of 
grading-restricted vertex superalgebras
and M\"{o}bius vertex superalgebras, respectively. We always work in the 
general setting of vertex superalgebras instead of just vertex algebras because
there is a canonical involution of a vertex superalgebra and
twisted modules associated to this involution are particularly important.
But we would like to emphasize that the results in this series of 
papers are new even for vertex algebras, not 
merely generalizations to the super case of some results on vertex algebras. 

The present paper is organized as follows: In the next section, we give a review of the construction 
in \cite{H-const-twisted-mod}. In Section 3, 
we give a  linearly independent
set of generators of the lower-bounded generalized $g$-twisted $V$-module 
$\widehat{M}^{[g]}_{B}$ constructed in Section 5 of \cite{H-const-twisted-mod}. 
In Section 4, we obtain the consequence
that the twisted Zhu's algebra or the twisted zero-mode algebra is not $0$. The main result
of this section is the existence 
of irreducible  lower-bounded generalized $g$-twisted $V$-modules. In Section 5, we prove that 
the weak commutativity for twisted fields in the assumption of 
Theorem 4.3 in \cite{H-const-twisted-mod} is in fact a consequence of the other assumptions.
In Section 6, we give several spanning sets of  lower-bounded generalized $g$-twisted $V$-modules
and discuss relations among elements of these spanning sets. We prove in Section 7 the existence of 
an irreducible grading-restricted generalized $g$-twisted $V$-module (an irreducible ordinary 
twisted module if the automorphism is of finite order) when $V$ is a M\"{o}bius vertex algebra 
and satisfies some conditions. 
In Section 8, we prove that a lower-bounded generalized $V^{g}$-module 
with an action of $g$   can be extended to a
lower-bounded generalized $g$-twisted $V$-module.

\renewcommand{\theequation}{\thesection.\arabic{equation}}
\renewcommand{\thethm}{\thesection.\arabic{thm}}
\setcounter{equation}{0}
\setcounter{thm}{0}
\section{Review of the construction of lower-bounded generalized $g$-twisted 
$V$-modules}

Since the results in this paper depend heavily on the construction of lower-bounded generalized $g$-twisted
$V$-modules given in \cite{H-const-twisted-mod}, we give a review of 
this construction in this section. 

We fix a grading-restricted vertex superalgebra $V$ and an automorphism $g$ of $V$ in this paper. 
Then $V=\coprod_{\alpha\in P_{V}}V^{[\alpha]}$,
where $V^{[\alpha]}$ is the generalized eigenspace for $g$ with eigenvalue $e^{2\pi i\alpha}$
and $P_{V}$ is the set $\{\alpha\in \C\;|\; \Re(\alpha)\in [0, 1), 
e^{2 \pi i\alpha} \;\text{is an eigenvalue of}\;g\}$. 
Then $g=e^{2\pi i \mathcal{L}_{g}}=e^{2\pi i (\mathcal{S}_{g}+\mathcal{N}_{g})}$,
where $\mathcal{L}_{g}$ is an operator on $V$ and $\mathcal{S}_{g}$ and $\mathcal{N}_{g}$
are the semisimple and nilpotent parts, respectively, of $\mathcal{L}_{g}$.
We assume that $V$ is generated by  
$\phi^{i}(x)=Y_{V}(\phi^{i}_{-1}\one, x)$ for $i\in I$, where
$\phi^{i}(x)$ or $\phi^{i}_{-1}\one$ for $i\in I$ are 
homogeneous with respect to weights and $\Z_{2}$-fermion numbers and
where 
$\phi^{i}_{-1}$ is the constant term of $\phi^{i}(x)$ and 
$\phi^{i}_{-1}\one=\lim_{x\to 0}\phi^{i}(x)\one$ (see \cite{H-2-const} for more details.)
For $i\in I$, $\phi^{i}_{-1}\one$ 
is a generalized eigenvector of $g$ with eigenvalue
$e^{2\pi i\alpha_{i}}$. 
We also assume that for $i\in I$, either  $\mathcal{N}_{g}\phi^{i}_{-1}\one=0$ 
or there exists 
$\mathcal{N}_{g}(i)\in I$ such that $\mathcal{N}_{g}\phi^{i}_{-1}\one
=\phi^{\mathcal{N}_{g}(i)}_{-1}\one$.

To formulate the construction theorem in \cite{H-const-twisted-mod}, we 
need some data and properties satisfied by these data. 

\renewcommand{\labelenumi}{$($\alph{enumi}$)$}
\begin{data}\label{data}
{\rm 
\begin{enumerate}

\item\label{data-a} Let 
$$W = \coprod_{n \in \C, s\in \Z_{2}, [\alpha]\in \C/\Z} W_{[n]}^{s; [\alpha]}
=\coprod_{n \in \C, s\in \Z_{2}, \alpha\in P_{W}} W_{[n]}^{s; [\alpha]}$$
 be a ${\C}\times \Z_{2}\times \C/\Z$-graded
vector space such that $W_{[n]}=\coprod_{s\in \Z_{2},  
\alpha\in P_{W}}W^{s; [\alpha]}_{[n]}=0$ 
when the real part of $n$ is sufficiently negative, where $P_{W}$ is the subset of 
the set $\{\alpha\in \C\;|\;\Re(\alpha)\in [0, 1)\}$ such that $W_{[n]}^{s; [\alpha]}\ne 0$
for $\alpha\in P_{W}$.

\item\label{data-b} Let 
\begin{align*}
\phi^{i}_{W}: W&\to x^{-\alpha^{i}}W((x))[\log  x]\nn
w&\mapsto \phi^{i}_{W}(x)w=\sum_{k\in \N}\sum_{n\in \alpha^{i}+\Z}
(\phi^{i}_{W})_{n, k}wx^{-n-1}
(\log x)^{k}
\end{align*}
for $i\in I$ be a set of linear maps called the {\it generating twisted field maps}.
Since $\phi^{i}_{W}(x)w\in x^{-\alpha^{i}}W((x))[\log  x]$, 
we must have $(\phi^{i}_{W})_{n, k}w=0$ when $n-\alpha^{i}$ is sufficiently 
negative and $k$ is sufficiently large. These linear maps correspond to
multivalued analytic maps with the preferred 
branch $\phi_{W}^{i; 0}$ and labeled branches $\phi_{W}^{i; p}$ for 
$p\in \Z$ from $\C^{\times}$
to $\hom(W, \overline{W})$. 

\item\label{data-c} Let 
\begin{align*}
\psi_{W}^{a}: V&\to \sum_{\alpha\in P_{V}}x^{-\alpha}W((x))[\log x]\nn
v&\mapsto \phi^{a}_{W}(x)v=\sum_{k\in \N}\sum_{\alpha\in P_{V},\; n\in \alpha+\Z}
(\psi^{a}_{W})_{n, k}vx^{-n-1}
(\log x)^{k}
\end{align*} 
for $a\in A$ be a set of linear maps called the {\it generator twist field maps}
such that $\phi^{a}_{W}(x)v\in x^{-\alpha}W((x))[\log x]$ for $\alpha\in P_{V}$ and 
$v\in V^{[\alpha]}$. 
Since $\phi^{a}_{W}(x)v\in x^{-\alpha}W((x))[\log  x]$ for $v\in V^{[\alpha]}$, 
we must have $(\psi^{a}_{W})_{n, k}v=0$ when $n-\alpha$ is sufficiently 
negative and $k$ is sufficiently large.
These linear maps corresponds to
multivalued analytic maps with preferred branch $\psi^{a; 0}$ and labeled branches
$\psi^{a; p}$ for $p\in \Z$ from $\C^{\times}$
to $\hom(V, \overline{W})$. 

\item\label{data-d} Let $L_{W}(0)$ and $L_{W}(-1)$ be operators on $W$. 

\item\label{data-e} An action of $g$ on $W$, denoted still by $g$, and 
an operator, still denoted by $\mathcal{L}_{g}$ and its
semisimple and nilpotent parts, still denoted by $\mathcal{S}_{g}$ and $\mathcal{N}_{g}$, 
respectively, on $W$ such that $g=e^{2\pi i\mathcal{L}_{g}}=
e^{2\pi i(\mathcal{S}_{g}+\mathcal{N}_{g})}$ on $W$.

\end{enumerate}}
\end{data}

These data are assumed to satisfy the following properties:

\renewcommand{\labelenumi}{\arabic{enumi}.}
\begin{assum}\label{basic-properties}
The space $W$, the generating twisted field maps  
 $\phi^{i}_{W}$ for $i\in I$,  the generator twist field maps
$\psi_{W}^{a}$ for $a\in A$, the operators $L_{W}(0)$, $L_{W}(-1)$,
$g$,  $\mathcal{L}_{g}$, $\mathcal{S}_{g}$ and $\mathcal{N}_{g}$ on $W$
in Data \ref{data} have the following properties:

\begin{enumerate}

\item\label{property-1} 
There exist semisimple and nilpotent operators 
$L_{W}(0)_{S}$ and $L_{W}(0)_{N}$ on $W$ such that
$L_{W}(0)=L_{W}(0)_{S}+L_{W}(0)_{N}$. 
For $i\in I$, 
$[L_{W}(0), \phi_{W}^{i}(x)]=z\frac{d}{dx}\phi_{W}^{i}(x)+(\wt \phi^{i})\phi_{W}^{i}(x)$.
For $a\in A$, there exists $(\wt \psi_{W}^{a})\in \C$ and, when $L_{W}(0)_{N}\psi_{W}^{a}(x)\ne 0$,
there exists $L_{W}(0)_{N}(a)\in A$ such that
$L_{W}(0)\psi_{W}^{a}(x)-\psi^{a}(x)L_{V}(0)
=x\frac{d}{dx}\psi_{W}^{a}(x)+(\wt \psi_{W}^{a})\psi_{W}^{a}(x)+\psi_{W}^{L_{W}(0)_{N}(a)}(x)$,
where $\psi_{W}^{L_{W}(0)_{N}(a)}(x)=0$ when $L_{W}(0)_{N}\psi_{W}^{a}(x)= 0$.

\item\label{property-2} For $i\in I$,
$[L_{W}(-1), \phi_{W}^{i}(x)]=\frac{d}{dx}\phi_{W}^{i}(x)$
and for $a\in A$,
$L_{W}(-1)\psi_{W}^{a}(x)-\psi_{W}^{a}(x)L_{V}(-1)=\frac{d}{dx}\psi_{W}^{a}(x)$.

\item\label{property-3} For $a\in A$, $\psi_{W}^{a}(x)\one\in W[[x]]$ and its constant terms
$\lim_{x\to 0}\psi_{W}^{a}(x)\one$ 
is homogeneous with respect to weights,  $\Z_{2}$-fermion number and $g$-weights.

\item\label{property-4} The vector space $W$ is spanned by elements of the form
$(\phi_{W}^{i_{1}})_{n_{1}, l_{1}}\cdots (\phi_{W}^{i_{k}})_{n_{k}, l_{k}}(\psi_{W}^{a})_{n, l}
v$ for $i_{1}, \dots, i_{k}\in I$, $a\in \A$
and $n_{1}\in \alpha^{i_{1}}+\Z, \dots, n_{k}\in \alpha^{i_{k}}+\Z$, 
$n\in \C$, $l_{1}, \dots, l_{k}, l\in \N$, $v\in V$.

\item\label{property-5} (i) For $i\in I$, $g\phi_{W}^{i;p+1}(z)g^{-1}=\phi_{W}^{i; p}(z)$.
(ii) For $i\in I$, $\phi_{W}^{i}(x)=
x^{-\mathcal{N}_{g}}(\phi_{W}^{i})_{0}(x)x^{\mathcal{N}_{g}}$ and for $a\in A$, 
$\psi_{W}^{a}(x)=(\psi_{W}^{a})_{0}(x)x^{-\mathcal{N}_{g}}$
where 
$(\phi_{W}^{i})_{0}(x)$ and $(\psi_{W}^{a})_{0}(x)$ are the constant terms 
of $\phi_{W}^{i}(x)$ and $\psi_{W}^{a}(x)$, respectively,  viewed
as power series of $\log x$ (with coefficients being series in powers of $x$). 
(iii)  For $i\in I$, $e^{2\pi i\mathcal{S}_{g}}
\phi_{W}^{i}(z)e^{-2\pi i\mathcal{S}_{g}}=e^{2\pi \alpha^{i}}\phi_{W}^{i}(z)$ and 
$[\mathcal{N}_{g}, \phi_{W}^{i}(z)]=\phi_{W}^{\mathcal{N}_{g}(i)}(z).$
(iv) For $a\in A$, 
there exists $\alpha^{a}\in P_{W}$ such that $(\psi_{W}^{a})_{n, 0}\one$ for 
$n\in -\N-1$ are generalized eigenvectors of $g$ 
with eigenvalue $e^{2\pi i \alpha^{a}}$. 

\item\label{property-6}  For $i, j\in I$, there exists $M_{ij}\in \Z_{+}$ such that
$$(x_{1}-x_{2})^{M_{ij}}\phi_{W}^{i}(x_{1})\phi_{W}^{j}(x_{2})
=(x_{1}-x_{2})^{M_{ij}}(-1)^{|\phi^{i}||\phi^{j}|}\phi_{W}^{j}(x_{2})\phi_{W}^{i}(x_{1}).
$$

\item\label{property-7}  For $i\in I$ and $a\in A$, there exists $M_{ia}\in \Z_{+}$ such that
\begin{align*}
&(x_{1}-x_{2})^{\alpha_{i}+M_{ia}}
(x_{1}-x_{2})^{\mathcal{N}_{g}}\phi_{W}^{i}(x_{1})(x_{1}-x_{2})^{-\mathcal{N}_{g}}
\psi_{W}^{a}(x_{2})\nn
&\quad =(-x_{2}+x_{1})^{\alpha_{i}+M_{ia}}
 (-1)^{|\phi^{i}||\psi^{a}|}
\psi_{W}^{a}(x_{2})
(-x_{2}+x_1)^{\mathcal{N}_{g}} \phi^{i}(x_{1})
(-x_{2}+x_1)^{-\mathcal{N}_{g}}.
\end{align*}

\end{enumerate}
\end{assum}

We want to define a twisted vertex operator map
\begin{align*}
Y^{g}_{W}: V\otimes W &\to W\{x\}[\log x],\nn
u\otimes w&\mapsto Y^{g}_{W}(u, x)w.
\end{align*}
Such a map is equivalent to a multivalued analytic map (denoted using the same notation)
\begin{align*}
Y^{g}_{W}: \C^{\times}&\to \hom(V\otimes W, \overline{W}),\nn
z&\mapsto Y^{g}_{W}(\cdot, z)\cdot: u\otimes w\mapsto Y^{g}_{W}(u, z)w
\end{align*}
with labeled branches 
\begin{align*}
(Y^{g}_{W})^{p}: \C^{\times}&\to \hom(V\otimes W, \overline{W}),\nn
z&\mapsto (Y^{g}_{W})^{p}(\cdot, z)\cdot: u\otimes w\mapsto (Y^{g}_{W})^{p}(u, z)w
\end{align*}
for $p\in \Z$. 
For $w'\in W'$, $w\in W$, $i_{1}, \dots, i_{k}\in I$, $m_{1}, \dots, m_{k}\in \Z$,
we define $(Y^{g}_{W})^{p}$ by
\begin{align*}
\langle w', &(Y^{g}_{W})^{p}(\phi^{i_{1}}_{m_{1}}\cdots \phi^{i_{k}}_{m_{k}}\one, z)
w\rangle\nn
&=\res_{\xi_{1}=0}\cdots\res_{\xi_{k}=0}
\xi_{1}^{m_{1}}\cdots\xi_{k}^{m_{k}}
F^{p}(\langle w', \phi_{W}^{i_{1}}(\xi_{1}+z)\cdots \phi_{W}^{i_{k}}(\xi_{k}+z)
w\rangle),
\end{align*}
where $F^{p}(\langle w', \phi_{W}^{i_{1}}(\xi_{1}+z)\cdots \phi_{W}^{i_{k}}(\xi_{k}+z)
w\rangle)$ is the branch labeled by $p\in \Z$  of the multivalued function obtained from 
the analytic extension of the absolutely convergent sum of the series 
$$\langle w', \phi_{W}^{i_{1}}(\xi_{1}+z)\cdots \phi_{W}^{i_{k}}(\xi_{k}+z)
w\rangle.$$

Then the following construction theorem is proved in \cite{H-const-twisted-mod} (see Theorem 4.3
in \cite{H-const-twisted-mod}):

\begin{thm}\label{const-thm}
The pair $(W, Y^{g}_{W})$ is a lower-bounded generalized $g$-twisted
$V$-module generated 
by $(\psi_{W}^{a})_{n, k}v$ for $a\in A$, $n\in \alpha+\Z$,
$k\in \N$, $v\in V^{[\alpha]}$ and $\alpha\in P_{V}$.
Moreover, this is the unique lower-bounded generalized 
$g$-twisted $V$-module structure on $W$ generated by $(\psi_{W}^{a})_{n, k}v$ 
for $a\in A$, $n\in \alpha+\Z$,
$k\in \N$, $v\in V^{[\alpha]}$ and $\alpha\in P_{V}$
such that $Y_{W}(\phi^{i}_{-1}\one, z)=\phi_{W}^{i}(z)$
for $i\in I$.
\end{thm}

In Section 5 of \cite{H-const-twisted-mod}, a universal lower-bounded generalized $g$-twisted 
$V$-module $\widehat{M}_{B}^{[g]}$ is constructed explicitly by applying this theorem
to the data obtained from a suitable vector space  and a 
real number. We also need this construction in this paper.

Let $M$ be a $\Z_{2}$-graded vector space (graded by $\Z_{2}$-fermion numbers). 
Assume that $g$ acts on $M$ and there is an operator 
$L_{M}(0)$ on $M$. Assume that there exist operators $\mathcal{L}_{g}$,
$\mathcal{S}_{g}$, $\mathcal{N}_{g}$ such that on $M$, $g=e^{2\pi i\mathcal{L}_{g}}$
and $\mathcal{S}_{g}$ and $\mathcal{N}_{g}$ are the semisimple and nilpotent, respectively,
parts of $\mathcal{L}_{g}$. Assume also that $M$ is a direct sum of generalized 
eigenspaces for the operator $L_{M}(0)$ and $L_{M}(0)$ can be decomposed as 
the sum of its semisimple part $L_{M}(0)_{S}$ and nilpotent part $L_{M}(0)_{N}$. 
Moreover, assume that the real parts of the eigenvalues of 
$L_{M}(0)$ has a lower bound. 
Let $\{w^{a}\}_{a\in A}$ be a basis of $M$ consisting of vectors homogeneous 
in weights, $\Z_{2}$-fermion numbers and 
$g$-weights (eigenvalues of $g$) such that for $a\in A$, 
either $L_{M}(0)_{N}w^{a}=0$ or there exists $L_{M}(0)_{N}(a)\in A$ 
such that $L_{M}(0)_{N}w^{a}=w^{L_{M}(0)_{N}(a)}$. For simplicity, when 
$L_{M}(0)_{N}w^{a}=0$, we shall use $w^{L_{M}(0)_{N}(a)}$ to denote $0$. 
Then for $a\in A$, we always have $L_{M}(0)_{N}w^{a}=w^{L_{M}(0)_{N}(a)}$.
For $a\in A$, let $\alpha^{a}\in \C$ such that $\Re(\alpha^{a})
\in [0, 1)$ and $e^{2\pi i\alpha^{a}}$ is the eigenvalue of $g$ for the generalized 
eigenvector $w^{a}$.

Let 
$$\hat{V}_{\phi}^{[g]}
=\coprod_{i\in I, k\in \N}\C\mathcal{N}_{g}^{k}\phi^{i}_{-1}\one
\otimes t^{\alpha^{i}}\C[t, t^{-1}]\oplus \C L_{0}\oplus \C L_{-1},$$
where $L_{0}$ and $L_{-1}$ are fixed abstract basis elements of a vector space 
$\C L_{0}\oplus \C L_{-1}$.  Let $T(\hat{V}_{\phi}^{[g]})$ be the tensor algebra of
$\hat{V}_{\phi}^{[g]}$ and let 
$$\phi_{\hat{V}_{\phi}^{[g]}}^{i}(x)
=\sum_{n\in \alpha^{i}+\Z}((x^{-\mathcal{N}_{g}}
\phi^{i}_{-1}\one)\otimes t^{n})x^{-n-1}
\in x^{-\alpha^{i}}\hat{V}_{\phi}^{[g]}[[x, x^{-1}]][\log x]$$
for $i\in I$.

For $i, j\in I$, we can always find $M_{i, j}\in \Z_{+}$ such that 
$x_{0}^{M_{i, j}}Y_{V}(\phi^{i}_{-1}\one, x_0)\phi^{j}_{-1}\one$
is a power series in $x_{0}$. For each pair $i, j\in I$, we choose $M_{i, j}$ to be 
the smallest of such positive integers. 
Let $J(\hat{V}_{\phi}^{[g]})$ be the ideal of 
$T(\hat{V}_{\phi}^{[g]})$ generated by the coefficients
of the formal series
\begin{align*}
(x_{1}-x_{2})^{M_{ij}}\phi_{\hat{V}_{\phi}^{[g]}}^{i}(x_{1})
\phi_{\hat{V}_{\phi}^{[g]}}^{j}(x_{2})
&-(-1)^{|\phi^{i}||\phi^{j}|}(x_{1}-x_{2})^{M_{ij}}
\phi_{\hat{V}_{\phi}^{[g]}}^{j}(x_{2})
\phi_{\hat{V}_{\phi}^{[g]}}^{i}(x_{1}),\\
L_{0}\phi_{\hat{V}_{\phi}^{[g]}}^{i}(x)-\phi_{\hat{V}_{\phi}^{[g]}}^{i}(x)L_{0}
&-x\frac{d}{dx}\phi_{\hat{V}_{\phi}^{[g]}}^{i}(x)-(\wt \phi^{i})\phi_{\hat{V}_{\phi}^{[g]}}^{i}(x),\\
L_{-1}\phi_{\hat{V}_{\phi}^{[g]}}^{i}(x)-\phi_{\hat{V}_{\phi}^{[g]}}^{i}(x)L_{-1}
&-\frac{d}{dx}\phi_{\hat{V}_{\phi}^{[g]}}^{i}(x)
\end{align*}
for $i, j\in I$, where the tensor product symbol $\otimes$ is omitted.
 Let $U(\hat{V}_{\phi}^{[g]})
=T(\hat{V}_{\phi}^{[g]})/J(\hat{V}_{\phi}^{[g]})$. 

Let
$$\widetilde{M}^{[g]}=\coprod_{\alpha\in P_{V}}
U(\hat{V}_{\phi}^{[g]})\otimes (M\otimes 
t^{\alpha}\C[t, t^{-1}])\otimes V^{[\alpha]}.$$
For $i\in I$, 
we have the formal series of operators on $\widetilde{M}^{[g]}$
$$\phi_{\widetilde{M}^{[g]}}^{i}(x)
=\sum_{n\in \alpha^{i}+\Z}
((x^{-\mathcal{N}_{g}}\phi^{i}_{-1}\one)\otimes t^{n})x^{-n-1}.$$
For $a\in A$ and $v\in V^{[\alpha]}$, let
$$\psi_{\widetilde{M}^{[g]}}^{a}(x)v
=\sum_{n\in \alpha+\Z}
(w^{a}\otimes t^{n}) x^{-\mathcal{N}_{g}}v x^{-n-1}.$$

Let $B\in \R$ such that $B\le \Re(\wt w)$ for any generalized eigenvector $w\in M$  of
$L_{M}(0)$. 
Let $J_{B}(\widetilde{M}^{[g]})$ be the $U(\hat{V}_{\phi}^{[g]})$-submodule
of $\widetilde{M}^{[g]}$ generated by elements of the following forms:
(i) $(\psi_{\widetilde{M}^{[g]}}^{a})_{n, 0}\one$ for  $a\in A$,
and $n\not\in -\N-1$; (ii) (\ref{element-form}) for $i_{1}, \dots, i_{l}\in I$,  
$n_{1}\in \alpha^{i_{1}}+\Z,
\dots, n_{l}\in \alpha^{i_{l}}+\Z$, $0\le k_{1}\le K^{i_{1}}, \dots, 0\le k_{l}\le K^{i_{l}}$, 
$m=0, -1$, $a\in A$, $n \in \alpha+\Z$, $0\le k\le K$, $v\in V^{[\alpha]}$,
$\alpha\in P_{V}$ such that 
$$\Re(\wt \phi^{i_{1}}-n_{1}-1+\cdots +\phi^{i_{l}}-n_{l}-1-m+\wt w^{a} -n-1+\wt v)<B.$$
Consider the quotient  $U(\hat{V}_{\phi}^{[g]})$-module
 $\widetilde{M}^{[g]}/J_{B}(\widetilde{M}^{[g]})$. 

For $i\in I$ and $a\in A$, let $M_{i, a}\in \Z_{+}$ be the smallest of 
$m\in \Z$ such that $m>\wt \phi^{i}-1+\Re(\wt w^{a})-B-\Re(\alpha^{i})$. 
Let $J(\widetilde{M}^{[g]}/J_{B}(\widetilde{M}^{[g]}))$ 
be the $U(\hat{V}_{\phi}^{[g]})$-submodule
of $\widetilde{M}^{[g]}/J_{B}(\widetilde{M}^{[g]})$ generated by 
 the coefficients of the formal series 
\begin{align*}
(x_{1}&-x_{2})^{\alpha^{i}+M_{i, a}} (x_{1}-x_{2})^{\mathcal{N}_{g}}
\phi_{\widetilde{M}^{[g]}}^{i}(x_{1}) 
(x_{1}-x_{2})^{-\mathcal{N}_{g}}
\psi_{\widetilde{M}^{[g]}}^{a}(x_{2})v\nn
&-(-1)^{|u||w|} (-x_{2}+x_{1})^{\alpha^{i}+M_{i, a}}
\psi_{\widetilde{M}^{[g]}}^{a}(x_{2})
(-x_{2}+x_1)^{\mathcal{N}_{g}} \phi^{i}(x_{1}) 
(-x_{2}+x_1)^{-\mathcal{N}_{g}}v,
\end{align*}
for $i\in I$, $a\in A$ and $v\in V^{[\alpha]}$ and the coefficients of the formal series 
\begin{align*}
L_{\widetilde{M}^{[g]}}(0)\psi_{\widetilde{M}^{[g]}}^{a}(x)v-
\psi_{\widetilde{M}^{[g]}}^{a}(x)L_{V}(0)v
&-x\frac{d}{dx}\psi_{\widetilde{M}^{[g]}}^{a}(x)v-
(\wt w^{a})\psi_{\widetilde{M}^{[g]}}^{a}(x)v-\psi_{\widetilde{M}^{[g]}}^{L_{M}(0)_{N}a}(x)v,
\\
L_{\widetilde{M}^{[g]}}(-1)\psi_{\widetilde{M}^{[g]}}^{a}(x)v&-
\psi_{\widetilde{M}^{[g]}}^{a}(x)L_{V}(-1)v
-\frac{d}{dx}\psi_{\widetilde{M}^{[g]}}^{a}(x)v
\end{align*}
for $a\in A$ and $v\in V$.
Let
$$\widehat{M}^{[g]}_{B}=(\widetilde{M}^{[g]}_{\ell}/J_{B}(\widetilde{M}^{[g]}))
/J(\widetilde{M}^{[g]}/J_{B}(\widetilde{M}^{[g]})).$$
We shall use
$\phi_{\widehat{M}^{[g]}_{B}}^{i}(x)$, $\psi_{\widehat{M}^{[g]}_{B}}^{a}(x)$,
$L_{\widehat{M}^{[g]}_{B}}(0)$
and $L_{\widehat{M}^{[g]}_{B}}(-1)$ to denote the series of operators and the operators
on $\widehat{M}^{[g]}_{B}$ induced from the corresponding series of 
operators and operators  on $\widetilde{M}^{[g]}$. 

Using Theorem 4.3 in \cite{H-const-twisted-mod} (Theorem \ref{const-thm} above), 
the following result is proved in \cite{H-const-twisted-mod}
(see Theorem 5.1 in \cite{H-const-twisted-mod}): 

\begin{thm}
The twisted fields 
$\phi_{\widehat{M}^{[g]}_{B}}^{i}$ for $i\in I$
generate a twisted vertex operator map 
$$Y^{g}_{\widehat{M}^{[g]}_{B}}: V\otimes \widehat{M}^{[g]}_{B}\to \widehat{M}^{[g]}_{B}\{x\}[\log x]$$
such that $(\widehat{M}^{[g]}_{B}, Y^{g}_{\widehat{M}^{[g]}_{B}})$ is a 
lower-bounded generalized $g$-twisted $V$-module. Moreover, this is the unique generalized 
$g$-twisted $V$-module structure on $\widehat{M}^{[g]}_{B}$ generated 
by the coefficients of $(\psi_{\widehat{M}^{[g]}_{B}}^{a})(x)v$ for $a\in A$ and $v\in V$
such that $Y^{g}_{\widehat{M}^{[g]}_{B}}(\phi^{i}_{-1}\one, z)=\phi_{\widehat{M}^{[g]}_{B}}^{i}(z)$
for $i\in I$.
\end{thm}

The lower-bounded generalized $g$-twisted $V$-module 
$(\widehat{M}^{[g]}_{B}, Y^{g}_{\widehat{M}^{[g]}_{B}})$ has the following universal property
(Theorem 5.2 in \cite{H-const-twisted-mod}):

\begin{thm}\label{univ-prop}
Let $(W, Y^{g}_{W})$ be a lower-bounded generalized $g$-twisted $V$-module
and $M_{0}$ a $\Z_{2}$-graded subspace of $W$ invariant under the actions of 
$g$, $\mathcal{S}_{g}$, $\mathcal{N}_{g}$, $L_{W}(0)$, 
$L_{W}(0)_{S}$ and
$L_{W}(0)_{N}$. Let $B\in \R$ such that $W_{[n]}=0$ when $\Re(n)<B$. 
Assume that there is a linear map $f: M\to M_{0}$ 
preserving the $\Z_{2}$-fermion number grading and commuting with the actions of 
$g$, $\mathcal{S}_{g}$, $\mathcal{N}_{g}$, $L_{W}(0)$ $(L_{\widehat{M}^{[g]}_{B}}(0))$, 
$L_{W}(0)_{S}$ $(L_{\widehat{M}^{[g]}_{B}}(0)_{S})$ and
$L_{W}(0)_{N}$ $(L_{\widehat{M}^{[g]}_{B}}(0)_{N})$. Then there exists a unique module 
map $\tilde{f}: \widehat{M}^{[g]}_{B}\to W$ such that $\tilde{f}|_{M}=f$. 
If $f$ is surjective and $(W, Y^{g}_{W})$ is generated by the 
coefficients of $(Y^{g})_{WV}^{W}(w_{0}, x)v$ for $w_{0}\in M_{0}$ and $v\in V$, 
where $(Y^{g})_{WV}^{W}$ 
is the twist vertex operator map obtained from $Y_{W}^{g}$ (see \cite{H-twist-vo}), then 
$\tilde{f}$ is surjective. 
\end{thm}

An immediate consequence of Theorem \ref{univ-prop} is:

\begin{cor}\label{quotient}
Let $(W, Y^{g}_{W})$ be a lower-bounded generalized $g$-twisted $V$-module generated by the 
coefficients of $(Y^{g})_{WV}^{W}(w, x)v$ for $w\in M$, where $(Y^{g})_{WV}^{W}$ 
is the twist vertex operator map obtained from $Y_{W}^{g}$ (see \cite{H-twist-vo})
and $M$ is a $\Z_{2}$-graded subspace of $W$  invariant under the actions of 
$g$, $\mathcal{S}_{g}$, $\mathcal{N}_{g}$, $L_{W}(0)$, 
$L_{W}(0)_{S}$ and
$L_{W}(0)_{N}$. Let $B\in \R$ such that $W_{[n]}=0$ when $\Re(n)<B$. 
Then there is a generalized $g$-twisted $V$-submodule $J$ of 
$\widehat{M}^{[g]}_{B}$ such that $W$ is equivalent as a 
lower-bounded generalized $g$-twisted $V$-module to the quotient 
module $\widehat{M}^{[g]}_{B}/J$.
\end{cor}

\renewcommand{\theequation}{\thesection.\arabic{equation}}
\renewcommand{\thethm}{\thesection.\arabic{thm}}
\setcounter{equation}{0}
\setcounter{thm}{0}
\section{A  linearly independent set of generators}

 In this section, we first give a set of generators for the lower-bounded generalized $g$-twisted 
$V$-module given by Theorem 4.3 in \cite{H-const-twisted-mod} (Theorem \ref{const-thm} above). 
Then we prove that this set of generators is in fact linearly independent 
for the lower-bounded generalized $g$-twisted 
$V$-module  $\widehat{M}^{[g]}_{B}$
constructed in Section 5 of \cite{H-const-twisted-mod} (see the preceding section for a brief review). 

For a vertex operator superalgebra or a  vertex superalgebra $V$ with a 
conformal vector $\omega$, a $g$-twisted generalized $V$-module $W$
being generated by a subset $M$ of $W$ means that $W$ is spanned by elements of the form 
$(Y_{W})_{n_{1}, k_{1}}(v_{1})\cdots (Y_{W})_{n_{l}, k_{l}}(v_{l})w$
for $v_{1}, \dots, v_{l}\in V$, $n_{1}, \dots, n_{l}\in \C$, $k_{1}, \dots, k_{l}\in \N$
and $w\in M$. Using the associativity for $Y_{W}$, it is easy to see that 
in this case, $W$ is also spanned by
the coefficients of $Y_{W}(v, x)w$
for $v\in V$ and $w\in M$.
Since $\omega\in V$, $L_{W}(0)w$ 
and $L_{W}(-1)w$ for $w\in W$ are  coefficients of $Y_{W}(\omega, x)w$. 

But for a grading-restricted vertex (super)algebra $V$ which does not have a specified conformal 
vector but does have two operators  $L_{V}(0)$ and $L_{V}(-1)$, we have the the following 
definition:

\begin{defn}\label{generated}
{\rm We say that a lower-bounded generalized $g$-twisted generalized $V$-module $W$
is {\it generated by a subset $M$ of $W$} if $W$ is spanned by 
$$(Y^{g}_{W})_{n_{1}, k_{1}}(v_{1})\cdots (Y^{g}_{W})_{n_{l}, k_{l}}(v_{l})
L_{W}(-1)^{p}L_{W}(0)^{q}w$$
for $v_{1}, \dots, v_{l}\in V$, $n_{1}, \dots, n_{l}\in \C$, $k_{1}, \dots, k_{l}, p, q\in \N$
and $w\in M$, or equivalently, by
the coefficients of $Y^{g}_{W}(v, x)
L_{W}(-1)^{p}L_{W}(0)^{q}w$
for $v\in V$, $p, q\in \N$,
and $w\in M$.  }
\end{defn}

\begin{rema}
{\rm Note that if $L_{W}(0)^{q}w\in M$ for $q\in \N$, 
$W$ is in fact spanned by elements of the form 
$$(Y^{g}_{W})_{n_{1}, k_{1}}(v_{1})\cdots (Y^{g}_{W})_{n_{l}, k_{l}}(v_{l})
L_{W}(-1)^{p}w$$
 or the coefficients of 
$Y^{g}_{W}(v, x)L_{W}(-1)^{p}w$. }
\end{rema}

\begin{thm}\label{twisted-mod-generators}
Let $W$ be the lower-bounded $g$-twisted $V$-module given by 
Theorem 4.3 in \cite{H-const-twisted-mod} (for example, $\widehat{M}^{[g]}_{B}$
constructed in Section 5 of \cite{H-const-twisted-mod}, see Theorem \ref{const-thm} and Section 2). 
Then $W$ (in particular, $\widehat{M}^{[g]}_{B}$)
is spanned by elements of the form
\begin{equation}\label{element-form}
(\phi^{i_{1}}_{W})_{n_{1}, k_{1}}
\cdots (\phi^{i_{l}}_{W})_{n_{l}, k_{l}}
L_{W}(-1)^{k}(\psi_{W}^{a})_{-1, 0}\one
\end{equation}
for $i_{1}, \dots, i_{l}\in I$,  
$n_{1}\in \alpha^{i_{1}}+\Z,
\dots, n_{l}\in \alpha^{i_{l}}+\Z$, $0\le k_{1}\le K^{i_{1}}, \dots, 0\le k_{l}\le K^{i_{l}}$, 
$a\in A$, $k\in \N$. In particular, $W$
is generated by $(\psi^{a}_{W})_{-1, 0}\one$
for $a\in A$.
\end{thm}
\pf
Let $W_{0}$ be the space spanned by elements the form (\ref{element-form}).
Since $V$ is generated by $\phi^{i}(x)$ for $i\in I$,  
by the associativity of $Y_{W}^{g}$,
$W_{0}$ is the submodule of $W$ 
generated by $(\psi^{a}_{W})_{-1, 0}\one$
for $w\in M$ and $k\in \N$. We need to prove that $W=W_{0}$.

From Section 4 of \cite{H-const-twisted-mod} (see also Section 2), 
we know that $W$ is generated by elements of the forms
$(\psi^{a}_{W})_{n, k}v$
for  $a\in A$, $v\in V^{[\alpha]}$, $n\in 
\alpha+\Z$ and $k\in \N$. We need only prove that such elements are in $W_{0}$. 
Since $V$ is spanned by elements of the form
$\phi^{i_{1}}_{n_{1}}\cdots \phi^{i_{l}}_{n_{i_{l}}}\one$, we need only consider $v$ 
of this form. 

We use induction on $l$. By construction, we have 
$(\psi^{a}_{W})_{n, k}\one=0$ for $n\not\in -\N-1$ or $k\ne 0$.
By the commutator formula 
$$L_{W}(-1)\psi^{a}_{W}(x)
-\psi^{a}_{W}(x)L_{V}(-1)
=\frac{d}{dx}\psi^{a}_{W}(x),$$
we obtain 
$$(\psi^{a}_{W})_{-n-1, 0}\one
=\frac{1}{n!}L(-1)^{n}(\psi^{a}_{W})_{-1, 0}\one$$
for $n\in \N$. 
Thus when $l=0$, $(\psi^{a}_{W})_{n, k}\one=0$ 
for $n\not\in -\N-1$ or $k\ne 0$ and 
$$(\psi^{a}_{W})_{-n-1, 0}\one
=\frac{1}{n!}L(-1)^{n}(\psi^{a}_{W})_{-1, 0}\one\in W_{0}$$
for $a\in A$ and $n\in \N$.

Now assume that 
$(\psi^{a}_{W})_{n, k}v$
is in $W_{0}$ for $a\in A$, $n\in \C$ and $k\in \N$. We need to show that 
$(\psi^{a}_{W})_{n, k}
\phi^{i}_{m}v$ for $i\in I$, $m\in \Z$ and $v\in V$
is also in $W_{0}$. We use the induction on $m$.
Since $V$ is lower bounded, there exists $N\in -\N$ such that 
$V_{(n)}=0$ for $n< N$. Then for homogeneous $v\in V$, 
$\phi^{\mathcal{N}_{g}^{p}(i)}_{m}v=0$ 
for $p\in \N$ and $m> m_{i, v}=\wt \phi^{i}-1+\wt v-N$. 
By the generalized weak commutativity for $\psi_{W}^{a}$, we obtain
\begin{align}\label{twisted-mod-space-1}
(x_{1}&-x_{2})^{\alpha^{i}+M_{i, a}} (x_{1}-x_{2})^{\mathcal{N}_{g}}
\phi_{W}^{i}(x_{1}) 
(x_{1}-x_{2})^{-\mathcal{N}_{g}}
\psi_{W}^{a}(x_{2})v\nn
&=(-1)^{|\phi^{i}||w^{a}|} (-x_{2}+x_{1})^{\alpha^{i}+M_{i, a}}
\psi_{W}^{a}(x_{2})
(-x_{2}+x_1)^{\mathcal{N}_{g}} \phi^{i}(x_{1}) 
(-x_{2}+x_1)^{-\mathcal{N}_{g}}v\nn
&=(-1)^{|\phi^{i}||w^{a}|} (-x_{2}+x_{1})^{\alpha^{i}+M_{i, a}}
\psi_{W}^{a}(x_{2})\sum_{p\in \N}\frac{1}{p!}
(\log (-x_{2}+x_1))^{p}[\overbrace{\mathcal{N}_{g}, \cdots, [\mathcal{N}_{g}}^{p},
 \phi^{i}(x_{1}) ]\cdots]v\nn
&=\sum_{p\in \N}
\frac{(-1)^{|\phi^{i}||w^{a}|}}{p!}  (-x_{2}+x_{1})^{\alpha^{i}+M_{i, a}}
\psi_{W}^{a}(x_{2})
(\log (-x_{2}+x_1))^{p}
 \phi^{\mathcal{N}_{g}^{p}(i)}(x_{1}) v.
\end{align}
We use induction on the smallest $p\in \Z_{+}$ such that $\mathcal{N}_{g}^{p}(i)=0$.
In the case that $\mathcal{N}_{g}^{p}(i)=0$ for $p\in \Z_{+}$, 
the coefficient of $x_{1}^{-m_{i, v}-1}$ in the right-hand side of 
(\ref{twisted-mod-space-1}) is 
$$\frac{(-1)^{|\phi^{i}||w^{a}|}}{p!}  (-x_{2})^{\alpha^{i}+M_{i, a}}
\psi_{W}^{a}(x_{2})
 \phi^{i}_{m_{i, v}} v.$$
Since the coefficients of the left-hand side of (\ref{twisted-mod-space-1}) are in $W_{0}$,
by taking the coefficient of $x_{1}^{-m_{i, v}-1}$ in both sides of 
(\ref{twisted-mod-space-1}), we see that the coefficients of 
$\psi_{W}^{a}(x_{2})
 \phi^{i}_{m_{i, v}} v$ is in $W_{0}$. 
Assume that $(\psi^{a}_{W})_{n, k}
\phi^{i}_{m_{i, v}}v$ is in $W_{0}$ when 
$\mathcal{N}_{g}^{p}(i)=0$ for $p\ge q\in \Z_{+}$.
Then in the case that $\mathcal{N}_{g}^{p}(i)=0$ for $p\ge q+1\in \Z_{+}$,
the coefficient of $x_{1}^{-m_{i, v}-1}$ in the right-hand side of 
(\ref{twisted-mod-space-1}) is
\begin{align}\label{twisted-mod-space-2}
&\frac{(-1)^{|\phi^{i}||w^{a}|}}{p!}  (-x_{2})^{\alpha^{i}+M_{i, a}}
\psi_{W}^{a}(x_{2})
 \phi^{i}_{m_{i, v}} v\nn
&\quad +\sum_{p=1}^{q}
\frac{(-1)^{|\phi^{i}||w^{a}|}}{p!}  (-x_{2})^{\alpha^{i}+M_{i, a}}
\psi_{W}^{a}(x_{2})
(\log (-x_{2}))^{p}
 \phi^{\mathcal{N}_{g}^{p}(i)}_{m_{i, v}} v.
\end{align}
By assumption, the coefficients of the second term in (\ref{twisted-mod-space-2}) 
are in $W_{0}$. Since the coefficients of the left-hand side of 
(\ref{twisted-mod-space-1}) are also in $W_{0}$, the coefficients of the first term 
in (\ref{twisted-mod-space-2}) 
are in $W_{0}$. This finishes the proof that $(\psi^{a}_{W})_{n, k}
\phi^{i}_{m_{i, v}}v\in W_{0}$.

Now assume that $(\psi^{a}_{W})_{n, k}
\phi^{i}_{m}v\in W_{0}$ for $m_{i, v}\ge m>m_{0}$. We have to prove 
$(\psi^{a}_{W})_{n, k}
\phi^{i}_{m_{0}}v\in W_{0}$. Again we use induction on the smallest 
$p\in \Z_{+}$ such that $\mathcal{N}_{g}^{p}(i)=0$.
In the case that $\mathcal{N}_{g}^{p}(i)=0$ 
for $p\in \Z_{+}$,  the coefficient of $x_{1}^{-m_{0}-1}$
in the right-hand side of (\ref{twisted-mod-space-1}) is 
\begin{align}\label{twisted-mod-space-3}
(-1)&^{|\phi^{i}||w^{a}|}
(-x_{2})^{\alpha^{i}+M_{i, a}}
\psi_{W}^{a}(x_{2})
 \phi^{i}_{m_{0}} v\nn
&+\sum_{k\in \Z_{+}}
(-1)^{|\phi^{i}||w^{a}|} \binom{\alpha^{i}+M_{i, a}}{k}
(-x_{2})^{\alpha^{i}+M_{i, a}-k}
\psi_{W}^{a}(x_{2})
 \phi^{i}_{m_{0}+k} v.
\end{align}
By assumption, the coefficients of the second term in (\ref{twisted-mod-space-3})
is in $W_{0}$. Since the coefficients of the left-hand side of 
(\ref{twisted-mod-space-1}) are also in $W_{0}$, the coefficients of the first term 
in (\ref{twisted-mod-space-3}) 
are in $W_{0}$. Assume that $(\psi^{a}_{W})_{n, k}
\phi^{i}_{m_{0}}v$ is in $W_{0}$ when 
$\mathcal{N}_{g}^{p}(i)=0$ for $p\ge q\in \Z_{+}$.
Then in the case that $\mathcal{N}_{g}^{p}(i)=0$ for $p\ge q+1\in \Z_{+}$,
the coefficient of $x_{1}^{-m_{0}-1}$ in the right-hand side of 
(\ref{twisted-mod-space-1}) is
\begin{align}\label{twisted-mod-space-4}
&(-1)^{|\phi^{i}||w^{a}|}
 (-x_{2})^{\alpha^{i}+M_{i, a}}
\psi_{W}^{a}(x_{2})
 \phi^{i}_{m_{0}} v\nn
&\quad +\sum_{k\in \Z_{+}}
(-1)^{|\phi^{i}||w^{a}|} \binom{\alpha^{i}+M_{i, a}}{k}
 (-x_{2})^{\alpha^{i}+M_{i, a}+k}
\psi_{W}^{a}(x_{2})
 \phi^{i}_{m_{0}+k} v\nn
&\quad +\sum_{p=1}^{q}\sum_{k\in \N}\sum_{n\in \N}
\frac{(-1)^{|\phi^{i}||w^{a}|}}{p!} \binom{\alpha^{i}+M_{i, a}}{k}c_{pn}(x_{2}, \log x_{2})
 (-x_{2})^{\alpha^{i}+M_{i, a}+k}
\psi_{W}^{a}(x_{2})
 \phi^{\mathcal{N}_{g}^{p}(i)}_{m_{0}+k+n} v,
\end{align}
where $c_{pn}(x_{2}, \log x_{2})$ for $p, n\in \N$ are the coefficients of 
$(\log (-x_{2}+x_1))^{p}$ as a powers series in $x_{1}$, that is,
$$\sum_{n\in \N}c_{pn}(x_{2}, \log x_{2})x_{1}^{n}
=(\log (-x_{2}+x_1))^{p}.$$
By assumption, the coefficients of the second and third terms in (\ref{twisted-mod-space-4}) 
are in $W_{0}$. Since the coefficients of the left-hand side of 
(\ref{twisted-mod-space-1}) are also in $W_{0}$, the coefficients of the first term 
in (\ref{twisted-mod-space-4}) 
are in $W_{0}$. This finishes the proof that $(\psi^{a}_{W})_{n, k}
\phi^{i}_{m_{0}}v\in W_{0}$ and thus also finishes the proof of our theorem.
\epfv

We now prove that the elements $(\psi_{\widehat{M}^{[g]}_{B}}^{a})_{-k-1, 0}\one$ 
of $\widehat{M}^{[g]}_{B}$ for $a\in A$ and $k\in \N$ 
are in fact linearly independent. 

\begin{thm}\label{ind-psi}
The elements 
$$(\psi^{a}_{\widehat{M}^{[g]}_{B}})_{-k-1, 0}\one=
\frac{1}{k!}L_{\widehat{M}^{[g]}_{B}}(-1)^{k}
(\psi^{a}_{\widehat{M}^{[g]}_{B}})_{-1, 0}\one$$
for $k\in \N$ and $a\in A$ of $\widehat{M}^{[g]}_{B}$ are linearly independent. 
\end{thm}
\pf
Since $w^{a}$ for $a\in A$ is a basis of $M$, they are in particular linearly independent. 
So in 
$$\widetilde{M}^{[g]}=\coprod_{\alpha\in P_{V}}
U(\hat{V}_{\phi}^{[g]})\otimes (M\otimes 
t^{\alpha}\C[t, t^{-1}])\otimes V^{[\alpha]},$$
the elements 
$(\psi^{a}_{\widetilde{M}^{[g]}})_{-k-1, 0}\one=(w^{a}\otimes t^{-k-1})\otimes \one$ 
for $a\in A$ and $k\in \N$ are linearly independent. 

Next we prove that $(\psi^{a}_{\widetilde{M}^{[g]}})_{-k-1, 0}\one$
 for $k\in \N$ and $a\in A$ are linearly 
independent in $\widetilde{M}^{[g]}/J_{B}(\widetilde{M}^{[g]})$. Assume that
a finite linear combination $w$
of $(\psi^{a})_{\widetilde{M}^{[g]}})_{-k-1, 0}\one$ for
$a\in A$ and $k\in \N$ is in $J_{B}(\widetilde{M}^{[g]})$.
From the definition in Section 5 of \cite{H-const-twisted-mod} (see also Section 2),
$J_{B}(\widetilde{M}^{[g]})$ is the $U(\hat{V}_{\phi}^{[g]})$-submodule of 
$\widetilde{M}^{[g]}$ generated by elements of the forms 
$(\psi_{\widetilde{M}^{[g]}}^{a})_{n, 0}\one$ for  $a\in A$, $n\not\in -\N-1$ and
$$(\phi^{i_{1}}_{\widetilde{M}^{[g]}})_{n_{1}, k_{1}}
\cdots (\phi^{i_{l}}_{\widetilde{M}^{[g]}})_{n_{l}, k_{l}}L_{\widetilde{M}^{[g]}}(0)^{p}
L_{\widetilde{M}^{[g]}}(-1)^{q}
(\psi_{\widetilde{M}^{[g]}}^{a})_{n, k}v$$
for $i_{1}, \dots, i_{l}\in I$,  
$n_{1}\in \alpha^{i_{1}}+\Z,
\dots, n_{l}\in \alpha^{i_{l}}+\Z$, $0\le k_{1}\le K^{i_{1}}, \dots, 0\le k_{l}\le K^{i_{l}}$, 
$p, q\in \N$,
$a\in A$, $n \in \alpha+\Z$, $0\le k\le K$, $v\in V^{[\alpha]}$,
$\alpha\in P_{V}$, where $K\in \N$ satisfying $\mathcal{N}_{g}^{K+1}v=0$,  
such that 
$$\Re(\wt \phi^{i_{1}}-n_{1}-1+\cdots +\phi^{i_{l}}-n_{l}-1-m+\wt w^{a} -n-1+\wt v)<B.$$
Since for $a\in A$ and $k\in \N$, the weight of 
$(\psi^{a}_{\widetilde{M}^{[g]}})_{-k-1, 0}\one$ is 
$\wt w^{a}+k\ge B$, $w$
can only be in the $U(\hat{V}_{\phi}^{[g]})$-submodule of 
$\widetilde{M}^{[g]}$ generated by elements of the form
$(\psi_{\widetilde{M}^{[g]}}^{a})_{n, 0}\one$ for  $a\in A$, $n\not\in -\N-1$.
But the intersection of the 
$U(\hat{V}_{\phi}^{[g]})$-submodules of $\widetilde{M}^{[g]}$ generated
by $(\psi^{a}_{\widetilde{M}^{[g]}})_{-k-1, 0}\one$ for $a\in A$ and $k\in \N$ and generated by 
$(\psi_{\widetilde{M}^{[g]}}^{a})_{n, 0}\one$ for $a\in A$ and $n\not\in -\N-1$ 
are $0$. Thus $w$ as an element of this intersection must be $0$. Since we have proved 
that $(\psi^{a}_{\widetilde{M}^{[g]}})_{-k-1, 0}\one$ for $a\in A$ and $k\in \N$ are linearly 
independent in $\widetilde{M}^{[g]}$, the coefficients of $w$ expressed as a linear combination of 
these linearly independent elements must be $0$. 

Finally we prove that $(\psi^{a}_{\widetilde{M}^{[g]}})_{-k-1, 0}\one$
 for $a\in A$ and $k\in \N$ are indeed linearly independent in $\widehat{M}^{[g]}_{B}$.
By definition, 
$$\widehat{M}^{[g]}_{B}=(\widetilde{M}^{[g]}_{\ell}/J_{B}(\widetilde{M}^{[g]}))
/J(\widetilde{M}^{[g]}/J_{B}(\widetilde{M}^{[g]})),$$
where 
$J(\widetilde{M}^{[g]}/J_{B}(\widetilde{M}^{[g]}))$ 
is the $U(\hat{V}_{\phi}^{[g]})$-submodule
of $\widetilde{M}^{[g]}/J_{B}(\widetilde{M}^{[g]})$ generated by 
 the coefficients of the formal series  
\begin{align}
(x_{1}-x_{2})^{\alpha^{i}+M_{i, a}} (x_{1}-x_{2})^{\mathcal{N}_{g}}
\phi_{\widetilde{M}^{[g]}}^{i}&(x_{1}) 
(x_{1}-x_{2})^{-\mathcal{N}_{g}}
\psi_{\widetilde{M}^{[g]}}^{a}(x_{2})v\nn
-(-1)^{|u||w|} (-x_{2}+x_{1})^{\alpha^{i}+M_{i, a}}
&\psi_{\widetilde{M}^{[g]}}^{a}(x_{2})
(-x_{2}+x_1)^{\mathcal{N}_{g}} \phi^{i}(x_{1}) 
(-x_{2}+x_1)^{-\mathcal{N}_{g}}v,\label{phi-psi-element}\\
L_{\widetilde{M}^{[g]}}(0)\psi_{\widetilde{M}^{[g]}}^{a}(x)v-
\psi_{\widetilde{M}^{[g]}}^{a}(x)L_{V}(0)v
&-x\frac{d}{dx}\psi_{\widetilde{M}^{[g]}}^{a}(x)v-
(\wt w^{a})\psi_{\widetilde{M}^{[g]}}^{a}(x)v
-\psi_{\widetilde{M}^{[g]}}^{L_{M}(0)_{N}(a)}(x)v,\label{psi-L(0)-commutator}\\
L_{\widetilde{M}^{[g]}}(-1)\psi_{\widetilde{M}^{[g]}}^{a}(x)v&-
\psi_{\widetilde{M}^{[g]}}^{a}(x)L_{V}(-1)v
-\frac{d}{dx}\psi_{\widetilde{M}^{[g]}}^{a}(x)v\label{psi-L(-1)-commutator}
\end{align}
for $i\in I$, $a\in A$ and $v\in V^{[\alpha]}$, $\alpha\in P(V)$.
Assume that a linear combination $w$ of $(\psi^{a}_{\widetilde{M}^{[g]}})_{-k-1, 0}\one
\in \widetilde{M}^{[g]}/J_{B}(\widetilde{M}^{[g]})$ for $a\in A$ and $k\in \N$
is in $J(\widetilde{M}^{[g]}/J_{B}(\widetilde{M}^{[g]}))$.
Note that (\ref{phi-psi-element}) contains 
$\phi_{\widetilde{M}^{[g]}}^{i}(x_{1})$ but not $L_{\widetilde{M}^{[g]}}(0)$
and $L_{\widetilde{M}^{[g]}}(-1)$,  (\ref{psi-L(0)-commutator})  
contains $L_{\widetilde{M}^{[g]}}(0)$ but not $\phi_{\widetilde{M}^{[g]}}^{i}(x_{1})$ 
and $L_{\widetilde{M}^{[g]}}(-1)$, and (\ref{psi-L(-1)-commutator})
contains $L_{\widetilde{M}^{[g]}}(-1)$ but not $\phi_{\widetilde{M}^{[g]}}^{i}(x_{1})$ 
and $L_{\widetilde{M}^{[g]}}(0)$. In particular, nonzero elements of 
the $U(\hat{V}_{\phi}^{[g]})$-submodule
$J(\widetilde{M}^{[g]}/J_{B}(\widetilde{M}^{[g]}))$
of $\widetilde{M}^{[g]}/J_{B}(\widetilde{M}^{[g]})$ generated by 
 the coefficients of (\ref{phi-psi-element}), 
(\ref{psi-L(0)-commutator}) and (\ref{psi-L(-1)-commutator})
must contain coefficients of 
$\phi_{\widetilde{M}^{[g]}}^{i}(x_{1})$, $L_{\widetilde{M}^{[g]}}(0)$
or $L_{\widetilde{M}^{[g]}}(-1)$. If $w\ne 0$ in 
$\widetilde{M}^{[g]}/J_{B}(\widetilde{M}^{[g]})$, then $w$ is  not in 
$J(\widetilde{M}^{[g]}/J_{B}(\widetilde{M}^{[g]}))$ since $w$  does not contain 
coefficients of $\phi_{\widetilde{M}^{[g]}}^{i}(x_{1})$, $L_{\widetilde{M}^{[g]}}(0)$
or $L_{\widetilde{M}^{[g]}}(-1)$. Contradiction. So
$w$ must be $0$ in $\widetilde{M}^{[g]}/J_{B}(\widetilde{M}^{[g]})$.
Since we have proved that $(\psi^{a}_{\widetilde{M}^{[g]}})_{-k-1, 0}\one$ 
for $a\in A$ and $k\in \N$ are linearly independent in $\widetilde{M}^{[g]}/J_{B}(\widetilde{M}^{[g]})$, 
the coefficients of $w$ expressed as a linear combination
of  $(\psi^{a}_{\widetilde{M}^{[g]}})_{-k-1, 0}\one$  for $a\in A$ and $k\in \N$ must be
$0$, proving the linear independence of $(\psi^{a}_{\widetilde{M}^{[g]}})_{-k-1, 0}\one$  
for $a\in A$ and $k\in \N$ in $\widehat{M}^{[g]}_{B}$.
\epfv

\renewcommand{\theequation}{\thesection.\arabic{equation}}
\renewcommand{\thethm}{\thesection.\arabic{thm}}
\setcounter{equation}{0}
\setcounter{thm}{0}
\section{Existence of irreducible lower-bounded generalized twisted modules}

One immediate consequence of Theorems \ref{ind-psi} is the nontrivial 
result that the
 lower-bounded generalized $g$-twisted $V$-module $\widehat{M}^{[g]}_{B}$ is not $0$. 
In this section, using this result, we prove the conjectures that the twisted Zhu's algebra or 
the twisted zero-mode algebra is not $0$ and that there exist  irreducible 
lower-bounded generalized twisted modules. 

First, we have the following result:

\begin{thm}\label{first-existence}
The lower-bounded generalized $g$-twisted $V$-module $\widehat{M}^{[g]}_{B}$ is not $0$. 
In particular, there exist nonzero lower-bounded generalized $g$-twisted $V$-modules.
\end{thm}
\pf
Since $(\psi^{a}_{\widehat{M}^{[g]}_{B}})_{-k-1, 0}\one$ for $a\in A$ and $k\in \N$ 
are linearly independent by Theorem \ref{ind-psi}, they are certainly not $0$. 
In particular, $\widehat{M}^{[g]}_{B}$ is not $0$.
\epfv

Before we show that the twisted Zhu's algebra $A_{g}(V)$ or equivalently
the twisted zero-mode algebra $Z_{g}(V)$ are not $0$,  
we need to discuss the relation between the 
notion of lower-bounded generalized $V$-module and 
the notion of $\overline{\C}_{+}$-graded
weak $g$-twisted $V$-module. To do this, we first need to 
give the correct notions of module map and equivalence between 
lower-bounded generalized $g$-twisted $V$-module. 

Given a lower-bounded generalized $g$-twisted $V$-module $W$, 
let $L^{C}_{W}(0)=L_{W}(0)+C$ for $C\in \C$. 
Then $W$ equipped with the same twisted vertex operator map $Y_{W}$ and the same 
operator $L_{W}(-1)$ but with $L_{W}(0)$ replaced by $L^{C}_{W}(0)=L_{W}(0)+C$
is also a lower-bounded generalized $g$-twisted $V$-module. Therefore we should 
view $W$ equipped with $L_{W}(0)$ and $W$ equipped with $L_{W}^{C}(0)$
as equivalent. In view of this, we have the following notions of module map and equivalence 
between  lower-bounded generalized $g$-twisted $V$-modules:

\begin{defn}\label{equivalence}
{\rm Let $W_{1}$ and 
$W_{2}$ be 
lower-bounded generalized $g$-twisted $V$-modules. In the case that 
$W_{2}$ is indecomposable, a {\it module map
from $W_{1}$ to
$W_{2}$} is a
linear map $f: W_{1}\to W_{2}$ preserving the $\Z_{2}$-fermion number grading,
commuting with the actions of $g$ and satisfying  
$f(Y^{g}_{W_{1}}(u, x)w)=Y^{g}_{W_{2}}(u, x)f(w)$,
$f(L_{W_{1}}(0)w)=L_{W_{2}}(0)f(w)+C_{f}f(w)$  and
$f(L_{W_{1}}(-1)w)=L_{W_{2}}(-1)f(w)$ for $w\in W_{1}$ and some $C_{f}\in \C$ 
independent of $w$. In the general case that 
$W_{2}$ is a direct sum of 
indecomposable submodules, a {\it module map
from $W_{1}$ to
$W_{2}$} is a
linear map $f: W_{1}\to W_{2}$ such that $f$ composed with the projections
to the indecomposable submodules are module maps to these submodules. 
A module map $f$ is said to be an {\it equivalence} if $f$ is invertible.}
\end{defn}

From the discussion and  definition above, we see that even 
for an indecomposable lower-bounded 
generalized $g$-twisted $V$-module $W$,  the weights of elements of $W$ 
can be shifted by any complex number and only the differences
of the eigenvalues of $L_{W}(0)$ are meaningful. 

Let $W$ be a lower-bounded generalized 
$g$-twisted $V$-module. We also need formulas for the commutators of 
$L_{W}(0)_{S}$ and $L_{W}(0)_{N}$ with $Y_{W}^{g}(u, x)$ for $u\in V$,
where $L_{W}(0)_{S}$ and $L_{W}(0)_{N}$ are the semisimple part and nilpotent part of 
$L_{W}(0)$, respectively.  
For $v\in V$, 
we have 
$$Y_{W}^{g}(v, x)=\sum_{k\in \N}\sum_{n\in \C}(Y_{W}^{g})_{n, k}(v)x^{-n-1}(\log x)^{k}.$$
Let 
$$Y_{W}^{g}(v, x, y)=\sum_{k\in \N}\sum_{n\in \C}(Y_{W}^{g})_{n, k}(v)x^{-n-1}y^{k}$$
for $v\in V$. 

\begin{prop}
Let $W$ be a lower-bounded generalized $g$-twisted $V$-module. Then we have
\begin{align}
[L_{W}(0)_{S}, Y_{W}^{g}(v, x)]&=x\frac{d}{dx}Y_{W}^{g}(v, x, y)\mbar_{y=\log x}
+Y_{W}(L_{V}(0)v, x),\label{L(0)_S-commutator}\\
[L_{W}(0)_{N}, Y_{W}^{g}(v, x)]&=\frac{d}{dy}Y_{W}^{g}(v, x, y)\mbar_{y=\log x}.
\label{L(0)_N-commutator}
\end{align}
\end{prop}
\pf
From the 
$L(0)$-commutator formula 
\begin{equation}\label{L(0)-commutator}
[L_{W}(0), Y_{W}^{g}(v, x)]=x\frac{d}{dx}Y_{W}(v, x)+Y_{W}(L_{V}(0)v, x),
\end{equation}
we obtain 
$$\left(L_{W}(0)-\wt v-\wt w-x\frac{d}{dx}\right)Y_{W}^{g}(v, x)w
=Y_{W}(v, x)(L_{W}(0)-\wt w)w$$
for homogeneous $v\in V$.
Taking the coefficient of $x^{-n-1}(\log x)^{k}$, we obtain
\begin{align*}
(L_{W}&(0)-(\wt v-n-1)-\wt w)(Y_{W}^{g})_{n, k}(v)w\nn
&=-(L_{W}(0)-\wt v-\wt w+k+1)(Y_{W}^{g})_{n, k+1}(v)w+(Y_{W}^{g})_{n, k}(v)(L_{W}(0)-\wt w)w.
\end{align*}
Thus we have 
\begin{align}\label{L(0)-semisimple}
(L_{W}&(0)-(\wt v-n-1)-\wt w)^{p}(Y_{W}^{g})_{n, k}(v)w\nn
&=\sum_{i=0}^{p}\binom{p}{i}\prod_{j=1}^{i}
(-(L_{W}(0)-\wt v-\wt w+k+j))(Y_{W}^{g})_{n, k+i}(v)(L_{W}(0)-\wt w)^{p-i}w
\end{align}
for $p\in \N$. Let $K_{1}$ be a positive integer such that $(L_{W}(0)-\wt w)^{K_{1}}w=0$. 
We also know that there exists $K_{2}\in \N$ such that 
$(Y_{W}^{g})_{n, k+i}(v)\tilde{w}=0$ when $i> K_{2}$ for all $\tilde{w}$ 
of the form $(L_{W}(0)-\wt w)^{p-i}w$. Let $p=K_{1}+K_{2}$ in (\ref{L(0)-semisimple}). Then 
the right-hand side of (\ref{L(0)-semisimple}) is $0$ and by  (\ref{L(0)-semisimple}), the left-hand side
is also $0$ . This shows that the weight of 
$(Y_{W}^{g})_{n, k}(v)w$ is $-(\wt v-n-1)-\wt w$. Thus we have 
$$L_{W}(0)_{S}(Y_{W}^{g})_{n, k}(v)w=(-(\wt v-n-1)-\wt w)(Y_{W}^{g})_{n, k}(v)w.$$
Since this formula holds for all $w$, $n$ and $k$, we obtain 
(\ref{L(0)_S-commutator}). Together with (\ref{L(0)-commutator}), 
(\ref{L(0)_S-commutator}) implies (\ref{L(0)_N-commutator}).
\epfv

The following result  is needed in Section 7:

\begin{prop}\label{semisimply}
Let $W$ be an irreducible 
lower-bounded generalized $g$-twisted $V$-module. Then
$L_{W}(0)$ 
acts on $W$ semisimply if and only if $g$ acts on $W$ semisimply.
\end{prop}
\pf
Assume that $g$ acts on $W$ semisimply. Then from the equivariance 
property, there cannot be terms containing the logarithm of the variable in 
the twisted vertex operators. In particular, the right-hand side of (\ref{L(0)_N-commutator})
is $0$. Thus by (\ref{L(0)_N-commutator}), $L_{W}(0)_{N}$ commutes with 
the twisted vertex operators. It is also clear that $L_{W}(0)_{N}$ commute with 
$L_{W}(0)$ and $L_{W}(-1)$.
The kernel of $L_{W}(0)_{N}$
 is not $0$ since $L_{W}(0)_{N}$ is nilpotent. Since $L_{W}(0)_{N}$ commutes with 
the twisted vertex operators, $L_{W}(0)$ and $L_{W}(-1)$, its kernel is
a submodule of $W$. But $W$ is irreducible and the kernel of 
 $L_{W}(0)_{N}$ is not $0$, the kernel of $L_{W}(0)_{N}$ must be $W$.
Thus $L_{W}(0)_{N}=0$.

Assume that $L_{W}(0)$ 
acts on $W$ semisimply. Then $L_{W}(0)_{N}=0$ and $L_{W}(0)=L_{W}(0)_{S}$. 
The $L(0)$-commutator formula (\ref{L(0)-commutator}) gives 
\begin{equation}\label{L(0)-conjugation}
(Y_{W}^{g})^{p}(u, z)
=e^{l_{p}(z)L_{W}(0)}(Y_{W}^{g})^{0}(e^{-l_{p}(z)L_{W}(0)}, 1)e^{-l_{p}(z)L_{W}(0)}
\end{equation}
for $u\in V$ and $p\in \Z$. Since the right-hand side of (\ref{L(0)-conjugation})
has no terms containing the logarithm of $z$, so does the left-hand side. 
Thus $Y_{W}^{g}(u, x)\in W\{x\}$. By the equivariance property, 
$g$ must acts on $W$ semisimply.
\epfv

We now give the relation between the 
notion of lower-bounded generalized $g$-twisted $V$-module (see \cite{H-twist-vo}) and 
the notion of $\overline{\C}_{+}$-graded
weak $g$-twisted $V$-module (see \cite{HY}). Note that in general 
a lower-bounded generalized $g$-twisted $V$-module with the given weight-grading 
might not be a $\overline{\C}_{+}$-graded
weak $g$-twisted $V$-module. But we can always shift the weight-grading. 

\begin{prop}\label{two-types-modules}
A lower-bounded generalized $g$-twisted $V$-module is equivalent to
a $\overline{\C}_{+}$-graded 
weak $g$-twisted $V$-module by changing the weight-grading by a real number. For  a $\overline{\C}_{+}$-graded
weak $g$-twisted $V$-module $W=\coprod_{n\in \overline{\C}_{+}}W_{n}$ with the twisted vertex operator map 
$Y_{W}^{g}$, let $L_{W}(0)_{S}$ be the operator on $W$ defined by 
$L_{W}(0)_{S}w=nw$ for $w\in W_{n}$ and
$L_{W}(0)_{N}$ and $L_{W}(-1)$ be operators on $W$ such that $L_{W}(0)_{N}$
is nilpotent and preserving the $\overline{\C}_{+}$-grading of $W$ and 
(\ref{L(0)_N-commutator}) and the $L(-1)$-commutator formula for $Y_{W}^{g}(u, x)$
holds for $u\in V$. Then $W$ equipped with the twisted vertex operator map 
$Y_{W}^{g}$ and the operators $L_{W}(0)=L_{W}(0)_{S}+L_{W}(0)_{N}$
and $L_{W}(-1)$ is a lower-bounded generalized $g$-twisted $V$-module.
\end{prop}
\pf
Let $W$ be a lower-bounded generalized $g$-twisted $V$-module
with the real parts of the weights of the homogeneous elements
being larger than or equal to $h$. We change 
$L_{W}(0)$ to $L_{W}^{-h}(0)=L_{W}(0)-h$.
Then $W$ equipped with $L_{W}^{-h}(0)$ is a 
lower-bounded generalized $g$-twisted $V$-module equivalent to $W$ 
but with the real parts of the weights of the homogeneous elements
being larger than or equal to $0$ so that the weights of its homogeneous subspaces are all in
$\overline{\C}_{+}$ (the closed right-half plane in the complex plane). 
In particular, it is a $\overline{\C}_{+}$-graded
weak $g$-twisted $V$-module. 

Let $W$ be a $\overline{\C}_{+}$-graded
weak $g$-twisted $V$-module. From the $\overline{\C}_{+}$-grading condition, 
(\ref{L(0)_S-commutator}) holds. Adding (\ref{L(0)_S-commutator}) and 
(\ref{L(0)_N-commutator}), we obtain the $L(0)$-commutator formula for $W$. 
Thus  $W$ equipped with  
$Y_{W}^{g}$, $L_{W}(0)=L_{W}(0)_{S}+L_{W}(0)_{N}$
and $L_{W}(-1)$ is a lower-bounded generalized $g$-twisted $V$-module.
\epfv

In \cite{DLM1}, Dong, Li and Mason generalized Zhu's algebra $A(V)$ associated to $V$ to 
a twisted Zhu's algebra $A_{g}(V)$ associated to $V$ and an automorphism $g$ of $V$ of finite order.
In \cite{HY}, Yang and the author introduced twisted zero-mode algebra $Z_{g}(V)$ associated to 
$V$ and an automorphism of $V$ not necessarily of finite order and also generalized 
the twisted Zhu's algebra $A_{g}(V)$ to the case that $g$ is not of finite order. These two associative algebras
are in fact isomorphic (see \cite{HY}). In the case that $g$ is of finite order, 
$A_{g}\ne 0$ is stated explicitly as a conjecture  in the beginning of Section 9 of the arxiv version 
of \cite{DLM1}. 
In the case that $V$ is $C_{2}$-cofinite and $g$ is of finite order, Dong, Li and Mason proved
this conjecture  in \cite{DLM2}. But in general, this conjecture has been open. The corollary below 
validates this conjecture:

\begin{cor}\label{z-g-not-0}
The twisted Zhu's algebra 
$A_{g}(V)$, or  equivalently, the twisted zero-mode algebra $Z_{g}(V)$, is not $0$. 
\end{cor}
\pf
Take $M$ in Section 5 of \cite{H-const-twisted-mod} (see also Section 2) to be $\C w^{a}$ for $a\in A=\{a\}$
such that $\wt w^{a}=B\in \overline{\C}_{+}$. 
Then $\widehat{M}^{[g]}_{B}$ is 
a $\overline{\C}_{+}$-graded weak $g$-twisted $V$-module. 
By Theorem \ref{first-existence},  $(\psi_{\widehat{M}^{[g]}_{B}}^{a})_{-1, 0}\one$  is not $0$. 
By definition, $(\psi_{\widehat{M}^{[g]}_{B}}^{a})_{-1, 0}\one\in \Omega_{g}(\widehat{M}^{[g]}_{B})$. 
So $\Omega_{g}(\widehat{M}^{[g]}_{B})$
is not $0$. Since $\Omega_{g}(\widehat{M}^{[g]}_{B})$ are nonzero $A_{g}(V)$- and $Z_{g}(V)$-modules,
$A_{g}(V)$ and  $Z_{g}(V)$ cannot be $0$.
\epfv

Next we prove the existence of 
irreducible lower-bounded generalized $g$-twisted $V$-module. 

\begin{thm}\label{second-existence}
Let $W$ be a lower-bounded generalized $g$-twisted $V$-module generated by 
a nonzero element $w$ (for example, $\widehat{M}^{[g]}_{B}$
when $M$ is a one dimensional space and $B$ is less than or equal to the real part of the weight 
of the elements of $M$). 
Then there exists a maximal submodule $J$ of $W$ such that 
$J$ does not contain $w$ and the quotient 
$W/J$ is irreducible.
\end{thm}
\pf
The set of submodules of $W$ not containing $w$ equipped with 
the relation of submodules
is a partially ordered set. For every  totally ordered subset of this set,  the union of 
all submodules in this subset is an upper bound in this set. Thus Zorn's lemma 
says that there must be a maximal submodule $J$ in this set of submodules. 

Since $w$ is not $0$ and $J$ does not contain $w$, $W/J$
is also not $0$. If there is a proper submodule of $W/J$, this submodule cannot contain 
$w+J$ since $w+J$ would generate $W/J$. Since $J$ also does not contain 
$w$, the inverse image of this proper submodule in $W/J$ under the projection from 
$W$ to $W/J$ must be a submodule of $W$
that does not contain $w$ but contains $J$. 
Since $J$ is maximal,  this inverse image must be in $J$
and thus must be equal to $J$. So the only proper submodule of  $W/J$ is $0$.
\epfv

\renewcommand{\theequation}{\thesection.\arabic{equation}}
\renewcommand{\thethm}{\thesection.\arabic{thm}}
\setcounter{equation}{0}
\setcounter{thm}{0}
\section{Weak commutativity for twisted fields as a consequence}

We prove in this section that 
the weak commutativity for twisted fields in the assumption of 
Theorem 4.3 in \cite{H-const-twisted-mod} (Theorem \ref{const-thm})
is in fact a consequence of the other assumptions.

We first give a conceptual motivation of this result.
As is mentioned in \cite{H-twist-vo},
the axioms for lower-bounded generalized $g$-twisted $V$-modules
are not independent. In fact, the commutativity for twisted vertex operators
follows from 
the associativity for twisted vertex operators, the commutativity for vertex operators for $V$ and 
other axioms. On the other hand, the associativity for twisted vertex operators
is in fact equivalent to the commtativity involving one twist vertex operator under
the assumption that the other axioms hold. 
Thus this commutativity and the commutativity for vertex operators for $V$
imply the commutativity for twisted vertex operators
when the other axioms hold. Also the commutativity involving one twist vertex operator 
follows from  the generalized weak commutativity involving 
one twist vertex operator. Moreover, the form of the correlation functions
and commutativity for twisted vertex operators imply
the weak commutativity for twisted vertex operators.
It is well known that the weak commutativity for vertex operators for $V$ implies the commutativity 
for vertex operators for $V$.   Thus we see that 
the generalized weak commutativity involving 
one twist vertex operator, the weak commutativity for vertex operators for $V$
and other axioms imply the weak commutativity for twisted 
vertex operators. 

Assumption 2.3 in 
\cite{H-const-twisted-mod} (Assumption \ref{basic-properties} in Section 2) is the condition for Theorem 4.3 
in \cite{H-const-twisted-mod} (Theorem \ref{const-thm} in Section 2) to be true.
Assumption 2.3 in 
\cite{H-const-twisted-mod} (Assumption \ref{basic-properties})
 includes
in particular both the weak commutativity for twisted generating fields
$\phi^{i}_{W}(x)$ (Property 6 in Assumption \ref{basic-properties}) and the generalized weak commutativity 
involving one twist field $\psi_{W}^{a}(x)$ (Property 7 in Assumption \ref{basic-properties}). 
Motivated by the discussion above,
it is reasonable to expect that Property 6 in Assumption \ref{basic-properties}
 is a consequence of Property 7, the weak commutativity 
for the generating fields $\phi^{i}(x)$ for $V$ and other properties 
 in Assumption \ref{basic-properties}. We now prove that this is indeed true. 

\begin{thm}\label{weak-comm-consq}
Assuming that  Data 2.2 in \cite{H-const-twisted-mod} (Data \ref{data} in Section 2) satisfying 
Properties 1--5 and 7 in Assumption 2.3 in \cite{H-const-twisted-mod} (Assumption 
\ref{basic-properties}). Then Property 6 in 
Assumption 2.3 in \cite{H-const-twisted-mod} (the 
weak commutativity of $\phi^{i}_{W}(x)$, see Assumption 
\ref{basic-properties}) holds. 
\end{thm}
\pf
For $i\in I$, there exists $K^{i}\in \Z_{+}$ such that 
$\mathcal{N}_{g}^{k}(i)=0$ for $k>K^{i}$. For $i, j\in I$, there exists
$N_{ij}\in \Z_{+}$ such that 
$$(x_{1}-x_{2})^{N_{ij}}\phi^{i}(x_{1})\phi^{j}(x_{2})
=(x_{1}-x_{2})^{N_{ij}}\phi^{i}(x_{2})\phi^{j}(x_{1}).$$
Let $M_{ij}$ be the maximum of the positive integers 
$N_{\mathcal{N}_{g}^{k}(i)\mathcal{N}_{g}^{l}(j)}$
for $k=0, \dots, K^{i}$ and $l=0, \dots, K^{j}$. Let $K^{ij}
=\max(K^{i}, K^{j})$. Then
for $i, j, i_{1}, \dots, i_{n}\in I$ and $a\in A$, by Property 7 in 
Assumption \ref{basic-properties} (the generalized weak commutativity 
for $\psi^{a}_{W}$), we have
\begin{align}\label{weak-comm-consq-1}
&(x-y)^{M_{ij}}(x-x_{n+1})^{\alpha^{i}+M_{ia}}
(y-x_{n+1})^{\alpha^{j}+M_{ja}}\cdot\nn
&\quad\quad \cdot
(x-x_{n+1})^{\mathcal{N}_{g}}\phi^{i}_{W}(x)(x-x_{n+1})^{-\mathcal{N}_{g}}
(y-x_{n+1})^{\mathcal{N}_{g}}\phi^{j}_{W}(y)(y-x_{n+1})^{-\mathcal{N}_{g}}\cdot\nn
&\quad\quad \cdot
\prod_{p=1}^{n}(x_{p}-x_{n+1})^{\alpha^{i_{p}}+M_{i_{p} a}}
(x_{p}-x_{n+1})^{\mathcal{N}_{g}}\phi^{i_{p}}_{W}(x_{p})(x_{p}-x_{n+1})^{-\mathcal{N}_{g}}
\psi^{a}_{W}(x_{n+1})\one\nn
&\; =(x-y)^{M_{ij}}(-x_{n+1}+x)^{\alpha^{i}+M_{ia}}
(-x_{n+1}+y)^{\alpha^{j}+M_{ja}}\cdot\nn
&\quad\quad \cdot
\psi^{a}_{W}(x_{n+1})(-x_{n+1}+x)^{\mathcal{N}_{g}}\phi^{i}_{W}(x)
(-x_{n+1}+x)^{-\mathcal{N}_{g}}
(-x_{n+1}+y)^{\mathcal{N}_{g}}\phi^{j}_{W}(y)(-x_{n+1}+y)^{-\mathcal{N}_{g}}\cdot\nn
&\quad\quad \cdot\prod_{p=1}^{n}(-x_{n+1}+x_{p})^{\alpha^{i_{p}}+M_{i_{p} a}}
(-x_{n+1}+x_{p})^{\mathcal{N}_{g}}\phi^{i_{p}}(x_{p})
(-x_{n+1}+x_{p})^{-\mathcal{N}_{g}}\one\nn
&\;=\sum_{k, l=0}^{K^{ij}}
\frac{1}{k!l!}(\log (-x_{n+1}+x))^{k}(\log (-x_{n+1}+y))^{l}
(-x_{n+1}+x)^{\alpha^{i}+M_{ia}}
(-x_{n+1}+y)^{\alpha^{j}+M_{ja}}\cdot\nn
&\quad\quad \cdot
\psi^{a}_{W}(x_{n+1})(x-y)^{M_{ij}}
\phi^{\mathcal{N}_{g}^{k}(i)}(x)\phi^{\mathcal{N}_{g}^{l}(j)}(y)\cdot\nn
&\quad\quad \cdot
\prod_{p=1}^{n}(-x_{n+1}+x_{p})^{\alpha^{i_{p}}+M_{i_{p} a}}
(-x_{n+1}+x_{p})^{\mathcal{N}_{g}}\phi^{i_{p}}(x_{p})
(-x_{n+1}+x_{p})^{-\mathcal{N}_{g}}\one\nn
&\;=\sum_{k, l=0}^{K^{ij}}
\frac{1}{k!l!}(\log (-x_{n+1}+x))^{k}(\log (-x_{n+1}+y))^{l}
(-x_{n+1}+x)^{\alpha^{i}+M_{ia}}
(-x_{n+1}+y)^{\alpha^{j}+M_{ja}}\cdot\nn
&\quad\quad \cdot
\psi^{a}_{W}(x_{n+1})(x-y)^{M_{ij}}\phi^{\mathcal{N}_{g}^{l}(j)}(y)
\phi^{\mathcal{N}_{g}^{k}(i)}(x)\cdot\nn
&\quad\quad \cdot
\prod_{p=1}^{n}(-x_{n+1}+x_{p})^{\alpha^{i_{p}}+M_{i_{p} a}}
(-x_{n+1}+x_{p})^{\mathcal{N}_{g}}\phi^{i_{p}}(x_{p})
(-x_{n+1}+x_{p})^{-\mathcal{N}_{g}}\one\nn
&\; =(x-y)^{M_{ij}}(-x_{n+1}+x)^{\alpha^{i}+M_{ia}}
(-x_{n+1}+y)^{\alpha^{j}+M_{ja}}\cdot\nn
&\quad\quad \cdot
\psi^{a}_{W}(x_{n+1})
(-x_{n+1}+y)^{\mathcal{N}_{g}}\phi^{j}_{W}(y)(-x_{n+1}+y)^{-\mathcal{N}_{g}}
(-x_{n+1}+x)^{\mathcal{N}_{g}}\phi^{i}_{W}(x)
(-x_{n+1}+x)^{-\mathcal{N}_{g}}\cdot\nn
&\quad\quad \cdot\prod_{p=1}^{n}(-x_{n+1}+x_{p})^{\alpha^{i_{p}}+M_{i_{p} a}}
(-x_{n+1}+x_{p})^{\mathcal{N}_{g}}\phi^{i_{p}}(x_{p})
(-x_{n+1}+x_{p})^{-\mathcal{N}_{g}}\one\nn
&\;=(x-y)^{M_{ij}}(x-x_{n+1})^{\alpha^{i}+M_{ia}}
(y-x_{n+1})^{\alpha^{j}+M_{ja}}\cdot\nn
&\quad\quad \cdot
(y-x_{n+1})^{\mathcal{N}_{g}}\phi^{j}_{W}(y)(y-x_{n+1})^{-\mathcal{N}_{g}}
(x-x_{n+1})^{\mathcal{N}_{g}}\phi^{i}_{W}(x)(x-x_{n+1})^{-\mathcal{N}_{g}}\cdot\nn
&\quad\quad \cdot
\prod_{p=1}^{n}(x_{p}-x_{n+1})^{\alpha^{i_{p}}+M_{i_{p} a}}
(x_{p}-x_{n+1})^{\mathcal{N}_{g}}\phi^{i_{p}}_{W}(x_{p})(x_{p}-x_{n+1})^{-\mathcal{N}_{g}}
\psi^{a}_{W}(x_{n+1})\one.
\end{align}
Dividing both sides of (\ref{weak-comm-consq-1}) by 
$$x^{\alpha^{i}+M_{ia}}y^{\alpha^{j}+M_{ja}}
\prod_{p=1}^{n}x_{p}^{\alpha^{i_{p}}+M_{i_{p} a}}$$
and then taking the constant term in $\log x$, $\log y$ and $\log x_{p}$ for $p=1, \dots, n$
in both sides,
we obtain 
\begin{align}\label{weak-comm-consq-2}
&(x-y)^{M_{ij}}\left(1-\frac{x_{n+1}}{x}\right)^{\alpha^{i}+M_{ia}}
\left(1-\frac{x_{n+1}}{y}\right)^{\alpha^{j}+M_{ja}}\cdot\nn
&\quad\quad \cdot
\left(1-\frac{x_{n+1}}{x}\right)^{\mathcal{N}_{g}}\phi^{i}_{W}(x)
\left(1-\frac{x_{n+1}}{x}\right)^{-\mathcal{N}_{g}}
\left(1-\frac{x_{n+1}}{y}\right)^{\mathcal{N}_{g}}\phi^{j}_{W}(y)
\left(1-\frac{x_{n+1}}{y}\right)^{-\mathcal{N}_{g}}\cdot\nn
&\quad\quad \cdot
\prod_{p=1}^{n}\left(1-\frac{x_{n+1}}{x_{p}}\right)^{\alpha^{i_{p}}+M_{i_{p} a}}
\left(1-\frac{x_{n+1}}{x_{p}}\right)^{\mathcal{N}_{g}}\phi^{i_{p}}_{W}(x_{p})
\left(1-\frac{x_{n+1}}{x_{p}}\right)^{-\mathcal{N}_{g}}
\psi^{a}_{W}(x_{n+1})\one\nn
&\quad =(x-y)^{M_{ij}}\left(1-\frac{x_{n+1}}{x}\right)^{\alpha^{i}+M_{ia}}
\left(1-\frac{x_{n+1}}{y}\right)^{\alpha^{j}+M_{ja}}\cdot\nn
&\quad\quad \cdot
\left(1-\frac{x_{n+1}}{y}\right)^{\mathcal{N}_{g}}\phi^{j}_{W}(y)
\left(1-\frac{x_{n+1}}{y}\right)^{-\mathcal{N}_{g}}
\left(1-\frac{x_{n+1}}{x}\right)^{\mathcal{N}_{g}}\phi^{i}_{W}(x)
\left(1-\frac{x_{n+1}}{x}\right)^{-\mathcal{N}_{g}}\cdot\nn
&\quad\quad \cdot
\prod_{p=1}^{n}\left(1-\frac{x_{n+1}}{x_{p}}\right)^{\alpha^{i_{p}}+M_{i_{p} a}}
\left(1-\frac{x_{n+1}}{x_{p}}\right)^{\mathcal{N}_{g}}\phi^{i_{p}}_{W}(x_{p})
\left(1-\frac{x_{n+1}}{x_{p}}\right)^{-\mathcal{N}_{g}}
\psi^{a}_{W}(x_{n+1})\one.
\end{align}
Note that by the definition of $M_{ij}$, 
(\ref{weak-comm-consq-2}) still holds when we replace 
$i$ and $j$ by $\mathcal{N}_{g}^{k}(i)$ and $\mathcal{N}_{g}^{l}(j)$,
respectively. 

We now  prove 
\begin{align}\label{weak-comm-consq-4}
(x-y)&^{M_{ij}}\phi^{i}_{W}(x)\phi^{j}_{W}(y)
\phi^{i_{1}}_{W}(x_{1})\cdots \phi^{i_{n}}_{W}(x_{n})
(\psi^{a}_{W})_{-q-1, 0}\one\nn
&=(x-y)^{M_{ij}}\phi^{j}_{W}(y)\phi^{i}_{W}(x)
\phi^{i_{1}}_{W}(x_{1})\cdots \phi^{i_{n}}_{W}(x_{n})
(\psi^{a}_{W})_{-q-1, 0}\one
\end{align}
for $q\in \N$ by using induction on $q$.
Note that $\psi^{a}_{W}(x_{n+1})\one$ is a power series in $x_{n+1}$. Taking the 
constant term in $x_{n+1}$ of both sides of (\ref{weak-comm-consq-2}), we obtain
\begin{align}\label{weak-comm-consq-3}
(x-y)&^{M_{ij}}\phi^{i}_{W}(x)\phi^{j}_{W}(y)
\phi^{i_{1}}_{W}(x_{1})\cdots \phi^{i_{n}}_{W}(x_{n})
(\psi^{a}_{W})_{-1, 0}\one\nn
&=(x-y)^{M_{ij}}\phi^{j}_{W}(y)\phi^{i}_{W}(x)
\phi^{i_{1}}_{W}(x_{1})\cdots \phi^{i_{n}}_{W}(x_{n})
(\psi^{a}_{W})_{-1, 0}\one.
\end{align}
This is (\ref{weak-comm-consq-4}) in the case $q=0$. Assume that 
(\ref{weak-comm-consq-4}) holds when $q< q_{0}$. Then in the case 
$q=q_{0}$, taking the coefficients of $x^{q_{0}}_{n+1}$ in both sides of 
(\ref{weak-comm-consq-2}), we obtain 
\begin{align*}
&(x-y)^{M_{ij}}\phi^{i}_{W}(x)\phi^{j}_{W}(y)
\phi^{i_{1}}_{W}(x_{1})\cdots \phi^{i_{n}}_{W}(x_{n})
(\psi^{a}_{W})_{-q_{0}-1, 0}\one\nn
&\;\quad +\res_{x_{n+1}}x_{n+1}^{-q_{0}-1}
(x-y)^{M_{ij}}\left(1-\frac{x_{n+1}}{x}\right)^{\alpha^{i}+M_{ia}}
\left(1-\frac{x_{n+1}}{y}\right)^{\alpha^{j}+M_{ja}}\cdot\nn
&\quad\quad \cdot
\left(1-\frac{x_{n+1}}{x}\right)^{\mathcal{N}_{g}}\phi^{i}_{W}(x)
\left(1-\frac{x_{n+1}}{x}\right)^{-\mathcal{N}_{g}}
\left(1-\frac{x_{n+1}}{y}\right)^{\mathcal{N}_{g}}\phi^{j}_{W}(y)
\left(1-\frac{x_{n+1}}{y}\right)^{-\mathcal{N}_{g}}\cdot\nn
&\quad\quad \cdot
\prod_{p=1}^{n}\left(1-\frac{x_{n+1}}{x_{p}}\right)^{\alpha^{i_{p}}+M_{i_{p} a}}
\left(1-\frac{x_{n+1}}{x_{p}}\right)^{\mathcal{N}_{g}}\phi^{i_{p}}_{W}(x_{p})
\left(1-\frac{x_{n+1}}{x_{p}}\right)^{-\mathcal{N}_{g}}\cdot\nn
\end{align*}
\begin{align}\label{weak-comm-consq-5}
&\quad\quad \cdot
\left(\sum_{q=0}^{k_{0}-1}(\psi^{a}_{W})_{-q-1, 0}\one x_{n+1}^{q}\right)\nn
&\; =(x-y)^{M_{ij}}\phi^{j}_{W}(y)\phi^{i}_{W}(x)
\phi^{i_{1}}_{W}(x_{1})\cdots \phi^{i_{n}}_{W}(x_{n})
(\psi^{a}_{W})_{-q_{0}-1, 0}\one\nn
&\;\quad +\res_{x_{n+1}}x_{n+1}^{-q_{0}-1}
(x-y)^{M_{ij}}\left(1-\frac{x_{n+1}}{x}\right)^{\alpha^{i}+M_{ia}}
\left(1-\frac{x_{n+1}}{y}\right)^{\alpha^{j}+M_{ja}}\cdot\nn
&\quad\quad \cdot
\left(1-\frac{x_{n+1}}{y}\right)^{\mathcal{N}_{g}}\phi^{j}_{W}(y)
\left(1-\frac{x_{n+1}}{y}\right)^{-\mathcal{N}_{g}}
\left(1-\frac{x_{n+1}}{x}\right)^{\mathcal{N}_{g}}\phi^{i}_{W}(x)
\left(1-\frac{x_{n+1}}{x}\right)^{-\mathcal{N}_{g}}\cdot\nn
&\quad\quad \cdot
\prod_{p=1}^{n}\left(1-\frac{x_{n+1}}{x_{p}}\right)^{\alpha^{i_{p}}+M_{i_{p} a}}
\left(1-\frac{x_{n+1}}{x_{p}}\right)^{\mathcal{N}_{g}}\phi^{i_{p}}_{W}(x_{p})
\left(1-\frac{x_{n+1}}{x_{p}}\right)^{-\mathcal{N}_{g}}\cdot\nn
&\quad\quad \cdot
\left(\sum_{q=0}^{q_{0}-1}(\psi^{a}_{W})_{-q-1, 0}\one x_{n+1}^{q}\right).
\end{align}
If the second terms in the two sides of 
(\ref{weak-comm-consq-5}) are equal, we obtain (\ref{weak-comm-consq-4}). 

We now use the induction assumption to prove that the second terms in the two sides of 
(\ref{weak-comm-consq-5}) are equal. The second terms in the left-hand
side and right-hand side of (\ref{weak-comm-consq-5}) are equal to 
\begin{align}\label{weak-comm-consq-6}
&\sum_{k, l=0}^{K^{ij}}\sum_{k_{1}=0}^{K^{i_{1}}}\cdots
\sum_{k_{n}=0}^{K^{i_{n}}}\sum_{q=0}^{q_{0}-1}
(x-y)^{M_{ij}}\phi^{\mathcal{N}_{g}^{k}(i)}_{W}(x)
\phi^{\mathcal{N}_{g}^{l}(j)}_{W}(y)\phi^{\mathcal{N}_{g}^{k_{1}}(i_{1})}_{W}(x_{1})\cdots
\phi^{\mathcal{N}_{g}^{k_{n}}(i_{n})}_{W}(x_{n})(\psi^{a}_{W})_{-q-1, 0}\one\cdot\nn
&\quad\quad \cdot
\res_{x_{n+1}}x_{n+1}^{q-k_{0}-1}
\left(1-\frac{x_{n+1}}{x}\right)^{\alpha^{i}+M_{ia}}
\left(1-\frac{x_{n+1}}{y}\right)^{\alpha^{j}+M_{ja}}\cdot\nn
&\quad\quad \cdot
\left(\log \left(1-\frac{x_{n+1}}{x}\right)\right)^{k}
\left(\log \left(1-\frac{x_{n+1}}{y}\right)\right)^{l}\cdot\nn
&\quad\quad \cdot
\prod_{p=1}^{n}\left(1-\frac{x_{n+1}}{x_{p}}\right)^{\alpha^{i_{p}}+M_{i_{p} a}}
\left(\log \left(1-\frac{x_{n+1}}{x_{p}}\right)\right)^{k_{p}}
\end{align}
and 
\begin{align}\label{weak-comm-consq-7}
&\sum_{k, l=0}^{K^{ij}}\sum_{k_{1}=0}^{K^{i_{1}}}\cdots
\sum_{k_{n}=0}^{K^{i_{n}}}\sum_{q=0}^{q_{0}-1}
(x-y)^{M_{ij}}\phi^{\mathcal{N}_{g}^{l}(j)}_{W}(y)
\phi^{\mathcal{N}_{g}^{k}(i)}_{W}(x)\phi^{\mathcal{N}_{g}^{k_{1}}(i_{1})}_{W}(x_{1})\cdots
\phi^{\mathcal{N}_{g}^{k_{n}}(i_{n})}_{W}(x_{n})(\psi^{a}_{W})_{-q-1, 0}\one\cdot\nn
&\quad\quad \cdot
\res_{x_{n+1}}x_{n+1}^{q-k_{0}-1}
\left(1-\frac{x_{n+1}}{x}\right)^{\alpha^{i}+M_{ia}}
\left(1-\frac{x_{n+1}}{y}\right)^{\alpha^{j}+M_{ja}}\cdot\nn
&\quad\quad \cdot
\left(\log \left(1-\frac{x_{n+1}}{x}\right)\right)^{k}
\left(\log \left(1-\frac{x_{n+1}}{y}\right)\right)^{l}\cdot\nn
&\quad\quad \cdot
\prod_{p=1}^{n}\left(1-\frac{x_{n+1}}{x_{p}}\right)^{\alpha^{i_{p}}+M_{i_{p} a}}
\left(\log \left(1-\frac{x_{n+1}}{x_{p}}\right)\right)^{k_{p}},
\end{align}
respectively. By the induction assumption, 
(\ref{weak-comm-consq-6}) and (\ref{weak-comm-consq-7}) are indeed equal.
This finishes our proof.
\epfv

\renewcommand{\theequation}{\thesection.\arabic{equation}}
\renewcommand{\thethm}{\thesection.\arabic{thm}}
\setcounter{equation}{0}
\setcounter{thm}{0}
\section{Spanning sets}

In this section, we give several spanning sets of $\widehat{M}^{[g]}_{B}$ and discuss
the relations among elements of these sets.

In Theorem \ref{twisted-mod-generators}, we obtain 
a spanning set of $\widehat{M}^{[g]}_{B}$ consisting 
of elements of 
the form (\ref{element-form}). This spanning set 
is certainly not linearly independent. But we still have the following result:

\begin{prop}\label{a-k-direct-sum}
For $a\in A$ and $k\in \N$, let  $\widehat{M}^{[g]}_{B; a, k}$ be the 
subspace of $\widehat{M}^{[g]}_{B}$ spanned by elements of the form 
(\ref{element-form}).
Then for $a_{0}\in A$ and $k_{0}\in \N$, the intersection of
$\widehat{M}^{[g]}_{B; a_{0}, k_{0}}$ with the subspace $\sum_{a\ne a_{0}, k\ne k_{0}}
\widehat{M}^{[g]}_{B; a, k}$ spanned  by elements
of the form (\ref{element-form}) for $a\ne a_{0}$ and $k\ne k_{0}$ is $0$.
In particular, $\widehat{M}^{[g]}_{B}$ is equal to the direct sum of $\widehat{M}^{[g]}_{B; a, k}$
for $a\in A$ and $k\in \N$.
\end{prop}
\pf
This proof is in fact a refinement of the proof of Theorem \ref{ind-psi}. 

In $\widetilde{M}^{[g]}$, for $a\in A$ and $k\in \N$, we have the subspace 
$\widetilde{M}^{[g]}_{a, k}$
 spanned by elements of the same form as (\ref{element-form}).
For $a_{0}\in A$ and $k_{0}\in \N$, the intersection of
$\widetilde{M}^{[g]}_{a_{0}, k_{0}}$ with the subspace $\sum_{a\ne a_{0}, k\ne k_{0}}
\widetilde{M}^{[g]}_{B; a, k}$ spanned  by elements
of the form (\ref{element-form}) in $\widetilde{M}^{[g]}$ for $a\ne a_{0}$ 
and $k\ne k_{0}$ is certainly $0$.

Next we prove that the intersection of the subspace 
$\widetilde{M}^{[g]}_{a_{0}, k_{0}}+J_{B}(\widetilde{M}^{[g]})$ of 
$\widetilde{M}^{[g]}/J_{B}(\widetilde{M}^{[g]})$
with the subspace $\sum_{a\ne a_{0}, k\ne k_{0}}
\widetilde{M}^{[g]}_{B; a, k}+J_{B}(\widetilde{M}^{[g]})$ 
is still $0$
in $\widetilde{M}^{[g]}/J_{B}(\widetilde{M}^{[g]})$. Let
$w+J_{B}(\widetilde{M}^{[g]})$ be in this intersection. We can take $w$ to be in 
$\widetilde{M}^{[g]}_{a_{0}, k_{0}}$. Then there exists 
$w_{1}\in J_{B}(\widetilde{M}^{[g]})$ such that  $w+w_{1}\in \sum_{a\ne a_{0}, k\ne k_{0}}
\widetilde{M}^{[g]}_{B; a, k}$. If $w$ is not in $J_{B}(\widetilde{M}^{[g]})$,
then $w+w_{1}$ is not $0$. Hence  $w+w_{1}\not\in \widetilde{M}^{[g]}_{a_{0}, k_{0}}$.
Since $w\in \widetilde{M}^{[g]}_{a_{0}, k_{0}}$, we see that 
$w_{1}\not\in \widetilde{M}^{[g]}_{a_{0}, k_{0}}$. On the other hand, 
from the definition in Section 5 of \cite{H-const-twisted-mod} (see Section 2),
$J_{B}(\widetilde{M}^{[g]})$ is the $U(\hat{V}_{\phi}^{[g]})$-submodule of 
$\widetilde{M}^{[g]}$ generated by elements of the form
$(\psi_{\widetilde{M}^{[g]}}^{a})_{n, 0}\one$ for  $a\in A$, $n\in \N$ and
$$(\phi^{i_{1}}_{\widetilde{M}^{[g]}})_{n_{1}, k_{1}}
\cdots (\phi^{i_{l}}_{\widetilde{M}^{[g]}})_{n_{l}, k_{l}}L_{\widetilde{M}^{[g]}}(0)^{p}
L_{\widetilde{M}^{[g]}}(-1)^{q}
(\psi_{\widetilde{M}^{[g]}}^{a})_{n, k}v$$
for $i_{1}, \dots, i_{l}\in I$,  
$n_{1}\in \alpha^{i_{1}}+\Z,
\dots, n_{l}\in \alpha^{i_{l}}+\Z$, $0\le k_{1}\le K^{i_{1}}, \dots, 0\le k_{l}\le K^{i_{l}}$, 
$p, q\in \N$,
$a\in A$, $n \in \alpha+\Z$, $0\le k\le K$, $v\in V^{[\alpha]}$,
$\alpha\in P_{V}$, where $K\in \N$ satisfying $\mathcal{N}_{g}^{K+1}v=0$,  
such that 
$$\Re(\wt \phi^{i_{1}}-n_{1}-1+\cdots +\phi^{i_{l}}-n_{l}-1-m+\wt w^{a} -n-1+\wt v)<B.$$
Then it is clear that we cannot add an element $w_{1}$ of the form above but not 
in $\widetilde{M}^{[g]}_{a_{0}, k_{0}}$ to $w$ to obtain an element of $\sum_{a\ne a_{0}, k\ne k_{0}}
\widetilde{M}^{[g]}_{B; a, k}$. Contradiction. Thus $w$ must be in $J_{B}(\widetilde{M}^{[g]})$
or equivalently, $w+J_{B}(\widetilde{M}^{[g]})$ is $0$ in $\widetilde{M}^{[g]}/J_{B}(\widetilde{M}^{[g]})$.
This proves that the intersection is $0$. 

Finally we prove that the intersection of
$\widehat{M}^{[g]}_{B; a_{0}, k_{0}}$ and $\sum_{a\ne a_{0}, k\ne k_{0}}
\widehat{M}^{[g]}_{B; a, k}$ is $0$. 
By definition, 
$$\widehat{M}^{[g]}_{B}=(\widetilde{M}^{[g]}_{\ell}/J_{B}(\widetilde{M}^{[g]}))
/J(\widetilde{M}^{[g]}/J_{B}(\widetilde{M}^{[g]})),$$
where 
$J(\widetilde{M}^{[g]}/J_{B}(\widetilde{M}^{[g]}))$ 
is the $U(\hat{V}_{\phi}^{[g]})$-submodule
of $\widetilde{M}^{[g]}/J_{B}(\widetilde{M}^{[g]})$ generated by 
 the coefficients of the formal series  (\ref{phi-psi-element}), 
(\ref{psi-L(0)-commutator}) and (\ref{psi-L(-1)-commutator})
for $i\in I$, $a\in A$ and $v\in V^{[\alpha]}$, $\alpha\in P(V)$.
For elements of the form (\ref{element-form}), the only relevant 
formal series are those with $v=\one$. But  (\ref{phi-psi-element}) has 
a term in which 
$\phi^{i}(x_{1})$ is to the right of $\psi_{W}^{a}(x_{2})$ so that it will not 
give any relations for elements of the form (\ref{element-form}). 
The series (\ref{psi-L(0)-commutator}) contain $L_{\widetilde{L}^{[g]}}(0)$
but elements of the form (\ref{element-form}) does not contain this operator.
So these series also do not give any relations for elements of the form (\ref{element-form}). 
Finally the series (\ref{psi-L(-1)-commutator}) in fact gives 
only the  relations 
$L_{\widehat{M}^{[g]}_{B}}(-1)(\psi_{\widehat{M}^{[g]}_{B}}^{a})_{-k-1, 0}\one
=(\psi_{\widehat{M}^{[g]}_{B}}^{a})_{-(k+1)-1, 0}\one$ for $k\in \N$ so that 
it also does not give relations for elements of the form (\ref{element-form}).  Hence adding elements
of $J(\widetilde{M}^{[g]}/J_{B}(\widetilde{M}^{[g]}))$  to 
elements of $\widetilde{M}^{[g]}_{a_{0}, k_{0}}+J_{B}(\widetilde{M}^{[g]})$
does not give elements of $\sum_{a\ne a_{0}, k\ne k_{0}}
\widetilde{M}^{[g]}_{B; a, k}+J_{B}(\widetilde{M}^{[g]})$. 
Thus the intersection of
$\widehat{M}^{[g]}_{B; a_{0}, k_{0}}$ and $\sum_{a\ne a_{0}, k\ne k_{0}}
\widehat{M}^{[g]}_{B; a, k}$ is $0$. 
\epfv

From Proposition \ref{a-k-direct-sum}, we see that all the 
relations among elements of the form (\ref{element-form}) are 
given by the relations among products of coefficients 
of $\phi_{\widehat{M}^{[g]}_{B}}^{i}(x)$ for $i\in I$. 
In particular, we can discuss such relations for fixed $a\in A$ and $k\in \N$. 

In the general setting of the present paper, since
the relations among $\phi^{i}(x)$ for $i\in I$ are not given explicitly, 
it is impossible to explicitly write down the relations among elements of 
this spanning set. But we can still give the existence of these relations corresponding to the 
relations among $\phi^{i}(x)$ for $i\in I$. 

Relations among the coefficients of $\phi^{i}(x)$ for $i\in I$ are of the form
\begin{equation}\label{relation-0}
\sum_{\nu=1}^{N}\lambda_{\nu}\phi^{i^{\nu}_{1}}_{m^{\nu}_{1}}
\cdots \phi^{i^{\nu}_{k}}_{m^{\nu}_{k}}\one=0
\end{equation}
for some $\lambda_{\nu}\in \C$ and $k\in \Z_{+}$, 
$i^{\nu}_{j}\in \tilde{I}$ (recall from Section 4 of \cite{H-const-twisted-mod} 
that $\tilde{I}=I\cup \{0\}$) and $m^{\nu}_{j}\in \Z$ for
$\nu=1, \dots, M$, $j=1, \dots, k$, such that 
$|\phi_{m_{1}^{\nu}}^{i_{1}^{\nu}}|+\cdots 
+|\phi_{m_{k}^{\nu}}^{i_{k}^{\nu}}|$ for 
$\nu=1, \dots, M$ are either all even or are all odd. 
In particular, for each $a\in A$ and $n\in \alpha^{i_{1}^{\nu}}+\cdots +
\alpha_{k}^{i_{k}^{\nu}}+\Z$, we have the relation 
\begin{equation}\label{relation}
\sum_{\nu=1}^{N}\lambda_{\nu}(\psi^{a}_{\widetilde{M}^{[g]}_{B}})_{n, 0}
\phi^{i^{\nu}_{1}}_{m^{\nu}_{1}}
\cdots \phi^{i^{\nu}_{k}}_{m^{\nu}_{k}}\one=0
\end{equation}
in $\widehat{M}^{[g]}_{B}$. 
By Theorem \ref{twisted-mod-generators},
the left-hand side of (\ref{relation}) can be rewritten as a linear combination of 
elements of the form (\ref{element-form}) such that 
$$\Re(\wt \phi^{i_{j}}-n_{j}-1+\cdots +\phi^{i_{l}}-n_{l}-1+\wt w^{a} +k)\ge B$$
for $j=1, \dots, l$. Thus the relation 
(\ref{relation}) can be rewritten as 
\begin{equation}\label{relation-1}
\sum_{\nu=1}^{N}\gamma_{\nu}
(\phi^{j^{\nu}_{1}}_{\widehat{M}^{[g]}_{B}})_{m^{\nu}_{1}, q^{\nu}_{1}}
\cdots (\phi^{j_{p}}_{\widehat{M}^{[g]}_{B}})_{m^{\nu}_{p}, q^{\nu}_{p}}
(\psi_{\widehat{M}^{[g]}_{B}}^{a})_{-q-1, 0}\one=0
\end{equation}
for $\gamma_{\nu}\in \C$, $j^{\nu}_{1}, \dots, j^{\nu}_{p}\in \tilde{I}$,  
$m^{\nu}_{1}\in \alpha^{j^{\nu}_{1}}+\Z,
\dots, m^{\nu}_{p}\in \alpha^{j^{\nu}_{p}}+\Z$, $0\le q^{\nu}_{1}
\le Q^{j^{\nu}_{1}}, \dots, 0\le q_{l}\le Q^{j^{\nu}_{p}}$, 
$q\in \N$ such that 
$$\Re(\wt \phi^{j^{\nu}_{r}}-m^{\nu}_{r}-1+\cdots 
+\phi^{j^{\nu}_{p}}-m^{\nu}_{p}-1+\wt w^{a} +q)\ge B$$ 
for $\nu=1, \dots, N$, $r=1, \dots, p$.
Since there are in general more than one way to rewrite 
the left-hand side of (\ref{relation}) as  linear combinations of 
elements of the form (\ref{element-form}), 
there are in general more than 
one relations of the form (\ref{relation-1}) corresponding to each relation of the form 
(\ref{relation-0}).

Thus we have 
the following result:

\begin{thm}\label{relation-thm}
The relations among the elements of the form (\ref{element-form})
are  all of the form (\ref{relation-1}),
for $j^{\nu}_{1}, \dots, j^{\nu}_{p}\in \tilde{I}$,  
$m^{\nu}_{1}\in \alpha^{j^{\nu}_{1}}+\Z,
\dots, m^{\nu}_{p}\in \alpha^{j^{\nu}_{p}}+\Z$, $0\le q^{\nu}_{1}
\le Q^{j^{\nu}_{1}}, \dots, 0\le q_{l}\le Q^{j^{\nu}_{p}}$, 
$q\in \N$ such that 
$$\Re(\wt \phi^{j^{\nu}_{r}}-m^{\nu}_{r}-1+\cdots 
+\phi^{j^{\nu}_{p}}-m^{\nu}_{p}-1+\wt w^{a} +q)\ge B$$ 
for $\nu=1, \dots, N$, $r=1, \dots, p$,  corresponding to all the relations of the form 
(\ref{relation-0}) in $V$.
\end{thm}
\pf
We have proved that all the relations of the form (\ref{relation-1}) are 
indeed satisfied by the elements of the form (\ref{element-form}).
There are apparently also another 
type of relations 
among elements of the form (\ref{element-form}) given by the 
coefficients of the weak commutativity for $\phi^{i}_{\widehat{M}^{[g]}_{B}}(x)$ for $i\in I$. 
But by Theorem \ref{weak-comm-consq},
these relations are obtained from the relations in $V$ given 
by the coefficients of the weak commutativity 
for $\phi^{i}(x)$ for $i\in I$. 
So the relations given by the coefficients of the weak commutativity 
for $\phi^{i}_{\widehat{M}^{[g]}_{B}}(x)$ for $i\in I$  are also of the form 
(\ref{relation-1}) corresponding to the relations in $V$ given by 
(\ref{relation-0}). 
Thus these are the only relations 
among such elements.
\epfv

Using the weak commutativity for $\phi^{i}_{\widehat{M}^{[g]}_{B}}(x)$ for $i\in I$, 
we now reduce the
spanning set in Proposition \ref{twisted-mod-generators}
to a smaller spanning set. 
We choose a total order on the index set $I$. 

\begin{lemma}\label{phi-order}
For $i_{1}, i_{2}\in I$ such that $i_{1}>i_{2}$, $n_{1}\in \alpha^{i_{1}}+\Z$,
$n_{2}\in \alpha^{i_{2}}+\Z$, $0\le k_{1}\le K^{i_{1}}$, $0\le k_{2}\le K^{i_{2}}$
and $w\in \widehat{M}^{[g]}_{B}$, 
the element
\begin{equation}\label{phi-order-1}
(\phi^{i_{1}}_{\widehat{M}^{[g]}_{B}})_{n_{1}, k_{1}}
(\phi^{i_{2}}_{\widehat{M}^{[g]}_{B}})_{n_{2}, k_{2}}w\in \widehat{M}^{[g]}_{B}
\end{equation}
can be written as a linear combination of elements of the form 
\begin{equation}\label{phi-order-2}
(\phi^{i_{2}}_{\widehat{M}^{[g]}_{B}})_{m_{2}, k_{2}}
(\phi^{i_{1}}_{\widehat{M}^{[g]}_{B}})_{m_{1}, k_{1}}w
\end{equation}
for $m_{1}\in \alpha^{i_{1}}+\Z$ and
$m_{2}\in \alpha^{i_{2}}+\Z$ such that $m_{1}-\alpha^{i_{1}}\le n_{1}-\alpha^{i_{1}}$,
$m_{2}-\alpha^{i_{2}}\ge n_{2}-\alpha^{i_{2}}$ and 
$m_{1}+m_{2}=n_{1}+n_{2}$.
\end{lemma}
\pf
We use induction on $n_{2}-\alpha^{i_{2}}$. Since $\widehat{M}^{[g]}_{B}$ is lower bounded,
there exists $N\in \Z$ such that $(\phi^{i_{2}}_{\widehat{M}^{[g]}_{B}})_{n_{2}, k_{2}}w=0$
for $n_{2}- \alpha^{i_{2}}>N$. 

In the case $n_{2}=\alpha^{i_{2}}+N$, 
$$(\phi_{\widehat{M}^{[g]}_{B}}^{i_{1}})_{n_{1}, k_{1}}
(\phi_{\widehat{M}^{[g]}_{B}}^{i_{2}})_{\alpha^{i_{2}}+N, k_{2}}w$$ 
is the coefficient of 
$x_{1}^{M_{i_{1}i_{2}}-n_{1}-1}(\log x_{1})^{k_{1}}
x_{2}^{-\alpha^{i_{2}}-N-1}(\log x_{2})^{k_{2}}$ in 
$$(x_{1}-x_{2})^{M_{i_{1}i_{2}}}\phi_{\widehat{M}^{[g]}_{B}}^{i_{1}}(x_{1})
\phi_{\widehat{M}^{[g]}_{B}}^{i_{2}}(x_{2})w.$$ 
By the weak commutativity for $\phi_{\widehat{M}^{[g]}_{B}}^{i_{1}}(x_{1})$
and $\phi_{\widehat{M}^{[g]}_{B}}^{i_{2}}(x_{2})$, we see that this coefficient
is equal to the coefficient of 
$x_{1}^{M_{i_{1}i_{2}}-n_{1}-1}(\log x_{1})^{k_{1}}x_{2}^{-\alpha^{i_{2}}-N-1}(\log x_{2})^{k_{2}}$ in 
$$(x_{1}-x_{2})^{M_{i_{1}i_{2}}}\phi_{\widehat{M}^{[g]}_{B}}^{i_{2}}(x_{2})
\phi_{\widehat{M}^{[g]}_{B}}^{i_{1}}(x_{1})w.$$ 
This coefficient is 
$$\sum_{p=0}^{M_{i_{1}i_{2}}}\binom{M_{i_{1}i_{2}}}{p}
(\phi_{\widehat{M}^{[g]}_{B}}^{i_{2}})_{\alpha^{i_{2}}+N+p, k_{2}}
(\phi_{\widehat{M}^{[g]}_{B}}^{i_{1}})_{n_{1}-p, k_{1}}w$$
and is indeed a linear combination of elements of the form (\ref{phi-order-2}). 

Now assume that for $n_{2}\in \alpha^{i_{2}}+\Z$
such that $n-\alpha^{i_{2}}<n_{2}-\alpha^{i_{2}}\le  N$, elements of the form 
(\ref{phi-order-1}) are linear combinations of elements of the form (\ref{phi-order-2}). 
In the case $n_{2}=n$, 
the coefficient of $x_{1}^{M_{i_{1}i_{2}}-n_{1}-1}(\log x_{1})^{k_{1}}
x_{2}^{-n-1}(\log x_{2})^{k_{2}}$ in 
$$(x_{1}-x_{2})^{M_{i_{1}i_{2}}}\phi_{\widehat{M}^{[g]}_{B}}^{i_{1}}(x_{1})
\phi_{\widehat{M}^{[g]}_{B}}^{i_{2}}(x_{2})w$$ 
is 
\begin{equation}\label{phi-order-3}
(\phi_{\widehat{M}^{[g]}_{B}}^{i_{1}})_{n_{1}, k_{1}}
(\phi_{\widehat{M}^{[g]}_{B}}^{i_{2}})_{n, k_{2}}w
+\sum_{p=1}^{M_{i_{1}i_{2}}}\binom{M_{i_{1}i_{2}}}{p}
(\phi_{\widehat{M}^{[g]}_{B}}^{i_{1}})_{n_{1}-p, k_{1}}
(\phi_{\widehat{M}^{[g]}_{B}}^{i_{2}})_{n+p, k_{2}}w.
\end{equation}
By the weak commutativity for $\phi_{\widehat{M}^{[g]}_{B}}^{i_{1}}(x_{1})$
and $\phi_{\widehat{M}^{[g]}_{B}}^{i_{2}}(x_{2})$, we see that 
(\ref{phi-order-3}) is equal to the coefficient of 
$x_{1}^{M_{i_{1}i_{2}}-n_{1}-1}(\log x_{1})^{k_{1}}
x_{2}^{-n-1}(\log x_{2})^{k_{2}}$ in 
$$(x_{1}-x_{2})^{M_{i_{1}i_{2}}}\phi_{\widehat{M}^{[g]}_{B}}^{i_{2}}(x_{2})
\phi_{\widehat{M}^{[g]}_{B}}^{i_{1}}(x_{1})w.$$ 
This coefficient is 
$$\sum_{p=0}^{M_{i_{1}i_{2}}}\binom{M_{i_{1}i_{2}}}{p}
(\phi_{\widehat{M}^{[g]}_{B}}^{i_{2}})_{n+p, k_{2}}
(\phi_{\widehat{M}^{[g]}_{B}}^{i_{1}})_{n_{1}-p, k_{1}}w$$
and is in fact a linear combination of elements of the form (\ref{phi-order-2}). 
Therefore (\ref{phi-order-3}) is also such a linear combination. 
But by  the induction assumption the second term in (\ref{phi-order-3}) is 
a linear combination of elements of the form (\ref{phi-order-2}). Thus 
the first term in (\ref{phi-order-3}) is 
a linear combination of elements of the form (\ref{phi-order-2}). 
By the principle 
of induction, the lemma is proved. 
\epfv

\begin{thm}\label{spanning-set}
The lower-bounded generalized $g$-twisted $V$-module  $\widehat{M}^{[g]}_{B}$
is spanned by elements of the form (\ref{element-form})
for $i_{1}, \dots, i_{l}\in I$ 
$n_{1}\in \alpha^{i_{1}}+\Z,
\dots, n_{l}\in \alpha^{i_{l}}+\Z$, $0\le k_{1}\le K^{i_{1}}, \dots, 0\le k_{l}\le K^{i_{l}}$, 
$a\in A$, $k\in \N$ such that  $i_{1}\le \cdots \le i_{l}$ and
$$\Re(\wt \phi^{i_{j}}-n_{j}-1+\cdots +\phi^{i_{l}}-n_{l}-1+\wt w^{a} +k)\ge B$$
for $j=1, \dots, l$ and $a\in  A$.  Moreover, the relations among these elements
are  of the form (\ref{relation-1}),
for $j^{\nu}_{1}, \dots, j^{\nu}_{p}\in \tilde{I}$,  
$m^{\nu}_{1}\in \alpha^{j^{\nu}_{1}}+\Z,
\dots, m^{\nu}_{p}\in \alpha^{j^{\nu}_{p}}+\Z$, $0\le q^{\nu}_{1}
\le Q^{j^{\nu}_{1}}, \dots, 0\le q_{l}\le Q^{j^{\nu}_{p}}$, 
$q\in \N$ such that $j^{\nu}_{1}\le \cdots \le j^{\nu}_{p}$ for $\nu=1, \dots, N$ and
$$\Re(\wt \phi^{j^{\nu}_{r}}-m^{\nu}_{r}-1+\cdots 
+\phi^{j^{\nu}_{p}}-m^{\nu}_{p}-1+\wt w^{a} +q)\ge B$$ 
for $\nu=1, \dots, N$, $r=1, \dots, p$,  corresponding to all the relations of the form 
(\ref{relation-0}) in $V$.
\end{thm}
\pf
This result follows immediately from Lemma \ref{phi-order} and the fact that
the real parts of the 
weights of the elements of $\widehat{M}^{[g]}_{B}$ are greater than or equal to $B$.
The proof of the second part is the same as the proof of Theorem \ref{relation-thm}.
\epfv

These spanning sets consist of elements where the operators 
$(\psi_{\widehat{M}^{[g]}_{B}})_{-k-1, 0}$ are to the right of all the other operators. 
We shall also give a spanning set where the operators $(\psi_{\widehat{M}^{[g]}_{B}})_{n, k}$
are to the left of all the other operators appearing in the elements. But we first need to 
give a spanning set of a different type.

\begin{prop}\label{span-1-v-o}
Let $W$ be a lower-bounded generalized $g$-twisted $V$-module generated 
by a set $S$ of homogeneous elements of $W$. Then $W$ is spanned by the coefficients of the formal
series of the form 
\begin{equation}\label{element-form-1}
Y_{W}^{g}(v, x)L_{W}(-1)^{k}w
\end{equation}
for $v\in V$, $k\in \N$ and $w\in S$. 
\end{prop}
\pf
Let $W_{0}$ be the subspace of $W$ spanned by elements of the form 
(\ref{element-form-1}). We prove that $W_{0}=W$. 

If $W_{0}$ is not equal to $W$, then there exists $\tilde{w}\in W$ but not in $W_{0}$.
Let $w'\in W'$ such that $\langle w', w_{0}\rangle=0$ for $w_{0}\in W_{0}$ and 
$\langle w', \tilde{w}\rangle=1$. Then for $v_{1}, \dots, v_{l}\in V$, $k\in \N$ and $w\in S$,
using the associativity, we have
\begin{align}\label{span-1-v-o-1}
F^{p}&(\langle w', Y_{W}^{g}(v_{1}, z_{1})\cdots Y_{W}^{g}(v_{l}, z_{l})
L_{W}(-1)^{k}w\rangle)\nn
&=F^{p}(\langle w', Y_{W}^{g}(Y_{V}(v_{1}, z_{1}-z_{l})\cdots Y_{V}(v_{l-1}, z_{l-1}-z_{l})
v_{l}, z_{l})L_{W}(-1)^{k}w\rangle).
\end{align}
But for fixed $z_{1}, \dots, z_{l}$ such that $|z_{l}|>|z_{1}-z_{l}|>\cdots
>|z_{l-1}-z_{l}|>0$,
$$Y_{W}^{g}(Y_{V}(v_{1}, z_{1}-z_{l})\cdots Y_{V}(v_{l-1}, z_{l-1}-z_{l})
v_{l}, z_{l})L_{W}(-1)^{k}w\in \overline{W}_{0}.$$
So by the definition of $w'$, the right-hand side of (\ref{span-1-v-o-1}) is $0$. 
By (\ref{span-1-v-o-1}), the left-hand side of (\ref{span-1-v-o-1})  
is also $0$. Thus 
\begin{equation}\label{span-1-v-o-2}
\langle w', Y_{W}^{g}(v_{1}, x_{1})\cdots Y_{W}^{g}(v_{l}, x_{l})L_{W}(-1)^{k}w\rangle=0
\end{equation}
for $v_{1}, \dots, v_{l}\in V$, $k\in \N$ and $w\in S$. Since 
$W$ is generated by $S$, (\ref{span-1-v-o-2}) for $v_{1}, \dots, v_{l}\in V$, $k\in \N$ and $w\in S$
implies that $w'=0$. But $\langle w', \tilde{w}\rangle=1$ so that $w'\ne 0$. Contradiction. 
Thus $W_{0}=W$.
\epfv

\begin{thm}\label{spanning-set-2}
The lower-bounded generalized $g$-twisted $V$ module  $\widehat{M}^{[g]}_{B}$
is spanned by the coefficients of the formal series
\begin{equation}\label{element-form-2}
Y_{\widehat{M}^{[g]}_{B}}^{g}(v, x)(\psi_{\widehat{M}^{[g]}_{B}}^{a})_{-k-1, 0}\one
=Y_{\widehat{M}^{[g]}_{B}}^{g}(v, x)L_{\widehat{M}^{[g]}_{B}}(-1)^{k}
(\psi_{W}^{a})_{-1, 0}\one
\end{equation}
for $v\in V$, $k\in \N$ and $a\in A$. In terms of the generating fields $\phi^{i}(x)$
of $V$, $\widehat{M}^{[g]}_{B}$
is spanned by the coefficients of the formal series
\begin{equation}\label{element-form-3}
Y_{\widehat{M}^{[g]}_{B}}^{g}(\phi^{i_{1}}_{n_{1}}\cdots
\phi^{i_{l}}_{n_{l}}\one, x)(\psi_{W}^{a})_{-k-1, 0}\one
\end{equation}
for $i_{1}, \dots, i_{l}\in I$, $n_{1}, \dots, n_{l}\in \Z$, $k\in \N$ and $a\in A$.
Moreover, the relations
among the coefficients of the formal series of the form (\ref{element-form-2}) 
are all induced from the relations among $v$'s in $V$ and the relations among the coefficients 
of the formal series (\ref{element-form-3})  are all induced from the relations among 
elements of the form $\phi^{i_{1}}_{n_{1}}\cdots
\phi^{i_{l}}_{n_{l}}\one$, that is, of the form (\ref{relation-0}).
\end{thm}
\pf
Since $\widehat{M}^{[g]}_{B}$ is generated by applying twisted vertex operators
to $(\psi_{W}^{a})_{-k-1, 0}\one$, by Proposition \ref{span-1-v-o}, we see that 
the coefficients of the formal series (\ref{element-form-2}) for $v\in V$, $k\in \N$ and $a\in A$
span $\widehat{M}^{[g]}_{B}$. Since $V$ is spanned by elements of the form 
$\phi^{i_{1}}_{n_{1}}\cdots
\phi^{i_{l}}_{n_{l}}\one$ for $i_{1}, \dots, i_{l}\in I$, $n_{1}, \dots, n_{l}\in \Z$,
$\widehat{M}^{[g]}_{B}$ is also spanned by elements of the form 
(\ref{element-form-3}). 

By Proposition \ref{a-k-direct-sum}, the relations 
among the coefficients 
of the formal series (\ref{element-form-3}) are all induced from the relations among 
elements of the form $\phi^{i_{1}}_{n_{1}}\cdots
\phi^{i_{l}}_{n_{l}}\one$, that is, of the form (\ref{relation-0}).
This implies that the relations
among the coefficients of the formal series of the form (\ref{element-form-2}) 
are all induced from the relations in $V$.
\epfv

\begin{thm}\label{spanning-set-3}
Let $W$ be a lower-bounded generalized $g$-twisted $V$-module. Let 
$\{\phi^{i}(x)\}_{i\in I}$ be a set of generating fields of $V$
and $\{\psi_{W}^{a}(x)\}_{a\in A}$ a set of generator twist fields of $W$.
Then $W$ is spanned by elements of the form 
$$(\psi_{W}^{a})_{n, k}\phi^{i_{1}}_{n_{1}}\cdots \phi^{i_{l}}_{n_{l}}\one$$
for $a\in A$, $i_{1}, \dots, i_{l}\in I$, $n\in \alpha^{i_{1}}+\cdots +\alpha^{i_{l}}+\Z$,
$k\in \N$,
$n_{1}, \dots, n_{l-1}\in \Z$, $n_{l}\in -\N-1$. 
\end{thm}
\pf
From the definition of the twist vertex operators in Section 4 of 
\cite{H-twist-vo}, for $a\in A$, $m\in \N$ and $v\in V$, we have 
\begin{align}\label{spanning-set-3-1}
Y^{g}_{W}&(v, x)(\psi_{W}^{a})_{-m-1, 0}\one\nn
&=Y^{g}_{W}(v, x)L_{W}(-1)^{m}(\psi_{W}^{a})_{-1, 0}\one\nn
&=(-1)^{|v||(\psi_{W}^{a})_{-1, 0}\one|}
e^{xL_{W}(-1)}(Y^{g})_{WV}^{W}(L_{W}(-1)^{m}(\psi_{W}^{a})_{-1, 0}\one, y)v
\mbar_{y^{n}=e^{-\pi ni}x^{n},\; \log y=\log x-\pi i}\nn
&=(-1)^{|v||(\psi_{W}^{a})_{-1, 0}\one|}
e^{xL_{W}(-1)}\frac{d^{m}}{dx^{m}}
(Y^{g})_{WV}^{W}((\psi_{W}^{a})_{-1, 0}\one, y)v
\mbar_{y^{n}=e^{-\pi ni}x^{n},\; \log y=\log x-\pi i}.
\end{align}
By Proposition \ref{span-1-v-o},
the  coefficients of the left-hand side of (\ref{spanning-set-3-1}) 
span $W$. So the coefficients of the right-hand side  of (\ref{spanning-set-3-1}) also span
$W$. But the coefficients of the right-hand side  of (\ref{spanning-set-3-1}) 
are linear combinations of $L_{W}(-1)^{m}(\psi_{W}^{a})_{n, k}u$ for $m\in \N$, 
$n\in \alpha+\Z$,
$k\in \N$, $u\in V^{[\alpha]}$ and $\alpha\in P_{V}$. From the $L(-1)$-commutator formula 
for the twist fields $\psi_{W}^{a}(x)$, we see that  
$L_{W}(-1)^{m}(\psi_{W}^{a})_{n, k}u$ is equal to a linear combination of elements of the 
form  $(\psi_{W}^{a})_{n', k'}L_{W}(-1)^{m'}u$ for $n'\in \alpha+\Z$,
$k' \in \N$ and $0\le m'\le m$. Thus elements 
of the form $(\psi_{W}^{a})_{n, k}u$ for $n\in \alpha+\Z$,
$k\in \N$, $u\in V^{[\alpha]}$ and $\alpha\in P_{V}$ span $W$. 
Since $u$ is spanned by elements of the form 
$\phi^{i_{1}}_{n_{1}}\cdots \phi^{i_{l}}_{n_{l}}\one$ for 
$i_{1}, \dots, i_{l}\in I$, $n_{1}, \dots, n_{l-1}\in \Z$, $n_{l}\in -\N-1$,
the result is proved.
\epfv

\renewcommand{\theequation}{\thesection.\arabic{equation}}
\renewcommand{\thethm}{\thesection.\arabic{thm}}
\setcounter{equation}{0}
\setcounter{thm}{0}
\section{Existence of irreducible grading-restricted generalized and ordinary
twisted modules}

In this section, as an application of Theorems \ref{second-existence} and \ref{spanning-set-2},
we prove the existence of irreducible grading-restricted generalized and ordinary 
$g$-twisted $V$-modules when $V$ is a M\"{o}bius vertex algebra and 
the twisted Zhu's algebra $A_{g}(V)$ or twisted zero-mode 
algebra $Z_{g}(V)$ is finite dimensional. In the case that $V$ is a vertex operator
algebra (meaning the existence of a conformal vector) and $g$ is of finite order, 
our result removed the assumption that $A_{g}(V)\ne 0$
in Theorem 9.1 in \cite{DLM1} and the assumption that $V$ is simple and 
$C_{2}$-cofinite in Theorem 9.1 in 
\cite{DLM2}. Our result is also much more general since $V$ does not have to have a conformal vector
(a M\"{o}bius structure or a compatible $\mathfrak{sl}_{2}$-module structure 
is enough) and $g$ can be of infinite order.

Note that Theorems \ref{first-existence} and \ref{second-existence} are about the 
existences of lower-bounded generalized 
$g$-twisted $V$-modules. We do not need any condition to prove these 
existences. We now want to discuss the existence of irreducible grading-restricted
generalized $g$-twisted $V$-modules. As is mentioned in the end of the introduction (Section 1), we 
call a grading-restricted generalized $g$-twisted $V$-module 
such that $L(0)$ acts semisimply an {\it ordinary $g$-twisted $V$-module}. We would also like to 
discuss the existence of irreducible ordinary
$g$-twisted $V$-modules.

To discuss these existences, we need to use contragredient modules and lowest 
weight subspaces. In particular, we
need to assume that our vertex superalgebra $V$ is a M\"{o}bius vertex superalgebra 
or a quasi-vertex superalgebra 
introduced in \cite{FHL} and \cite{HLZ} (see also \cite{H-twist-vo}) to
discuss the existence of lowest weights and lowest weight subspaces.
We first recall the definition of M\"{o}bius vertex superalgebra.

\begin{defn}
{\rm A {\it M\"{o}bius vertex superalgebra} is a grading-restricted vertex 
superalgebra $V$ equipped with 
an operator $L_{V}(1)$ satisfying
\begin{align*}
[L_{V}(1), L_{V}(0)]&=L_{V}(1),\\
[L_{V}(1), L_{V}(-1)]&=2L_{V}(0),\\
[L_{V}(1), Y_{V}(v, x)]&=Y_{V}(L_{V}(1)v, x)+2xY_{V}(L_{V}(0)v, x)
+x^{2}Y_{V}(L_{V}(-1)v, x)
\end{align*}
for $v\in V$. Let $V_{1}$ and $V_{2}$ be M\"{o}bius vertex superalgebras.
A {\it homomorphism $V_{1}$ to $V_{2}$} is a homomorphism $f$ from 
$V_{1}$ to $V_{2}$ when $V_{1}$ and $V_{2}$ are viewed as grading-restricted 
vertex superalgebras such that $f(L_{V_{1}}(1)v)=L_{V_{2}}(1)f(v)$. An {\it isomorphism
from $V_{1}$ to $V_{2}$} is an invertible homomorphism from $V_{1}$ to $V_{2}$.
An {\it automorphism of a M\"{o}bius vertex superalgebra $V$} is an isomorphism from $V$ 
to itself.}
\end{defn}

In this section, $V$ is a M\"{o}bius vertex superalgebra and $g$ is an automorphism of $V$.
In this case, a lower-bounded generalized $g$-twisted $V$-module
should also have an operator $L_{W}(1)$ and satisfies the corresponding properties.
Thus we need the following definition:

\begin{defn}\label{module-mobius}
{\rm Let $V$ be a M\"{o}bius vertex superalgebra. 
A {\it lower-bounded generalized $g$-twisted $V$-module} is a lower-bounded 
generalized $g$-twisted $V$-module for the underlying grading-restricted vertex 
superalgebra of $V$
equipped with an operator $L_{W}(1)$ on $W$ satisfying
\begin{align*}
[L_{W}(1), L_{W}(0)]&=L_{W}(1),\\
[L_{W}(1), L_{W}(-1)]&=2L_{W}(0),\\
[L_{W}(1), Y_{W}^{g}(v, x)]&=Y_{W}^{g}(L_{V}(1)v, x)+2xY_{W}^{g}(L_{V}(0)v, x)
+x^{2}Y_{W}^{g}(L_{V}(-1)u, x)
\end{align*}
for $v\in V$. Let $W_{1}$ and $W_{2}$ be lower-bounded generalized $g$-twisted $V$-module.
A {\it module map from $W_{1}$ to  $W_{2}$} is a module map $f$ from 
$W_{1}$ to $W_{2}$ when $W_{1}$ and $W_{2}$ are viewed as 
 lower-bounded generalized $g$-twisted $V$-modules for the underlying grading-restricted 
vertex superalgebra of $V$ such that $f(L_{W_{1}}(1)w)=L_{W_{2}}(1)f(w)$. An {\it equivalence
from $W_{1}$ to $W_{2}$} is an invertible module map from $W_{1}$ to $W_{2}$. }
\end{defn}

\begin{rema}
{\rm Note that for a lower-bounded generalized $g$-twisted $V$-module $W$,
because $L_{W}(0)$ appears in the right-hand side of the commutator formula 
between $L_{W}(1)$ and $L_{W}(-1)$, we cannot redefine $L_{W}(0)$ by adding 
a number to obtain a lower-bounded generalized $g$-twisted $V$-module with a different $L_{W}(0)$.
Thus in this case, the values, not just the differences, of the weights (the eigenvalues of 
$L_{W}(0)$) are meaningful. In particular, a module map $f: W_{1}\to W_{2}$
must satisfy $f(L_{W_{1}}(0)w)=L_{W_{2}}(0)f(w)$ for $w\in W_{1}$.}
\end{rema}

\begin{rema}
{\rm For a M\"{o}bius vertex superalgebra $V$, $L_{V}(-1)$, $L_{V}(0)$ and $L_{V}(1)$ 
give a structure of $\mathfrak{sl}_{2}$-module to $V$. Similarly, for a 
lower-bounded generalized $g$-twisted $V$-module, $L_{W}(-1)$, $L_{W}(0)$ and $L_{W}(1)$ 
give a structure of an $\mathfrak{sl}_{2}$-module to $W$.}
\end{rema}

\begin{rema}\label{second-existence-mobius}
{\rm Theorem \ref{second-existence} still holds for 
lower-bounded generalized $g$-twisted $V$-module $W$ in this case. In fact
the same proof works except that submodules means submodules invariant under
$L_{W}(1)$. }
\end{rema}

We also need to modify Definition \ref{generated} when $V$ is a 
M\"{o}bius vertex superalgebra. 

\begin{defn}\label{generated-mobius}
{\rm Let $V$ be a M\"{o}bius vertex superalgebra and $g$ an automorphism of $V$.
We say that a lower-bounded generalized $g$-twisted generalized $V$-module $W$
is {\it generated by a subset $M$ of $W$} if $W$ is spanned by 
$$(Y^{g}_{W})_{n_{1}, k_{1}}(v_{1})\cdots (Y^{g}_{W})_{n_{l}, k_{l}}(v_{l})
L_{W}(-1)^{p}L_{W}(0)^{q}L_{W}(1)^{r}w$$
for $v_{1}, \dots, v_{l}\in V$, $n_{1}, \dots, n_{l}\in \C$, $k_{1}, \dots, k_{l}, p, q, r\in \N$
and $w\in M$, or equivalently, by
the coefficients of $Y^{g}_{W}(v, x)
L_{W}(-1)^{p}L_{W}(0)^{q}L_{W}(1)^{r}w$
for $v\in V$, $p, q, r\in \N$,
and $w\in M$.  }
\end{defn}

Now we use the construction of lower-bounded generalized twisted 
modules for a grading-restricted vertex superalgebra
in \cite{H-const-twisted-mod} (see Section 2) to give a construction of 
lower-bounded generalized $g$-twisted 
modules for the M\"{o}bius vertex superalgebra $V$. We assume that the space $M$ in Section 5 of 
\cite{H-const-twisted-mod} (see Section 2) has an additional operator $L_{M}(1)$ of weight 
$-1$ and fermion number $0$ and commuting with $\mathcal{L}_{g}$,
$\mathcal{S}_{g}$ and $\mathcal{N}_{g}$. We also assume that there is a 
basis $\{w^{a}\}_{a\in A}$ 
of homogeneous basis of $M$ satisfying the condition given in Section 5 of 
\cite{H-const-twisted-mod} (that is, $L_{M}(0)_{N}w^{a}=w^{L_{M}(0)_{N}(a)}$
for $a\in A$, where either $L_{M}(0)_{N}(a)\in A$ or $w^{L_{M}(0)_{N}(a)}=0$, see Section 2)
and in addition satisfying the condition that $L_{M}(1)w^{a}=w^{L_{M}(1)(a)}$
for $a\in A$, where either $L_{M}(1)(a)\in A$ or $w^{L_{M}(1)(a)}=0$. 
Then we have the lower-bounded 
generalized $g$-twisted $V$-module $\widehat{M}^{[g]}_{B}$ when $V$ is viewed as a 
grading-restricted vertex superalgebra. We now define an operator $L_{\widehat{M}^{[g]}_{B}}(1)$
on $\widehat{M}^{[g]}_{B}$. 

By Theorem \ref{twisted-mod-generators}, $\widehat{M}^{[g]}_{B}$
is generated by $(\psi^{a}_{W})_{-1, 0}\one$
for $a\in A$. By Theorem \ref{spanning-set-2}, 
$\widehat{M}^{[g]}_{B}$ is spanned by the coefficients of the formal series 
of the form (\ref{element-form-2})
for $v\in V$, $k\in \N$ and $a\in A$. We define $L_{\widehat{M}^{[g]}_{B}}(1)$
by 
\begin{align*}
L_{\widehat{M}^{[g]}_{B}}(1)&Y_{\widehat{M}^{[g]}_{B}}^{g}(v, x)
L_{\widehat{M}^{[g]}_{B}}(-1)^{k}
(\psi^{a}_{\widehat{M}^{[g]}_{B}})_{-1, 0}\one\nn
&=Y_{\widehat{M}^{[g]}_{B}}^{g}((L_{V}(1)+2xL_{V}(0)+x^{2}L_{V}(-1))v, x)
L_{\widehat{M}^{[g]}_{B}}(-1)^{k}
(\psi^{a}_{\widehat{M}^{[g]}_{B}})_{-1, 0}\one\nn
&\quad +k(k-1)Y_{\widehat{M}^{[g]}_{B}}^{g}(v, x)
L_{\widehat{M}^{[g]}_{B}}(-1)^{k-1}(\psi^{a}_{\widehat{M}^{[g]}_{B}})_{-1, 0}\one\nn
&\quad +2k(\wt w^{a})Y_{\widehat{M}^{[g]}_{B}}^{g}(v, x)
L_{\widehat{M}^{[g]}_{B}}(-1)^{k-1}(\psi^{a}_{\widehat{M}^{[g]}_{B}})_{-1, 0}\one\nn
&\quad +Y_{\widehat{M}^{[g]}_{B}}^{g}(v, x)
L_{\widehat{M}^{[g]}_{B}}(-1)^{k}
(\psi^{L_{M}(1)(a)}_{\widehat{M}^{[g]}_{B}})_{-1, 0}\one
\end{align*}
for $v\in V$, $k\in \N$ and $a\in A$. By Theorem \ref{spanning-set-2}, 
 the only relations among elements of the form (\ref{element-form-2})
are given by the relations among $v$'s. Thus $L_{\widehat{M}^{[g]}_{B}}(1)$
is well defined.

\begin{thm}\label{2.5th-existence}
Let $V$ be a M\"{o}bius vertex superalgebra. 
Then $\widehat{M}^{[g]}_{B}$ equipped with 
$L_{\widehat{M}^{[g]}_{B}}(1)$ is a lower-bounded generalized $g$-twisted $V$-module.
Moreover, $\widehat{M}^{[g]}_{B}$ has the universal property stated as in 
Theorem 5.2 in \cite{H-const-twisted-mod} (Theorem \ref{univ-prop})
except that now $V$ is a M\"{o}bius vertex superalgebra.
Also every lower-bounded generalized $g$-twisted $V$-module generated by a $\Z_{2}$-graded subspace $M$ invariant 
under the actions of $g, \mathcal{S}_{g}, \mathcal{N}_{g}, L_{W}(0), L_{W}(0)_{S}, L_{W}(0)_{N}$
and bounded below in the real parts of the weights by $B\in \R$ is a quotient of $\widehat{M}^{[g]}_{B}$.
\end{thm}
\pf
From the definition, we have 
\begin{align*}
[L_{\widehat{M}^{[g]}_{B}}(1), L_{\widehat{M}^{[g]}_{B}}(0)]&=L_{\widehat{M}^{[g]}_{B}}(1),\\
[L_{\widehat{M}^{[g]}_{B}}(1), L_{\widehat{M}^{[g]}_{B}}(-1)]&=2L_{\widehat{M}^{[g]}_{B}}(0),\\
[L_{\widehat{M}^{[g]}_{B}}(1), Y_{\widehat{M}^{[g]}_{B}}^{g}(v, x)]
&=Y_{\widehat{M}^{[g]}_{B}}^{g}(L_{V}(1)v, x)+2xY_{\widehat{M}^{[g]}_{B}}^{g}(L_{V}(0)v, x)
+x^{2}Y_{\widehat{M}^{[g]}_{B}}^{g}(L_{V}(-1)v, x).
\end{align*}
By Definition \ref{module-mobius}, $\widehat{M}^{[g]}_{B}$ equipped with 
$L_{\widehat{M}^{[g]}_{B}}(1)$ is a lower-bounded generalized $g$-twisted $V$-module.

The universal properties and the statement that a lower-bounded generalized $g$-twisted $V$-module
is a  quotient follows immediately from the construction of 
$\widehat{M}^{[g]}_{B}$, Theorem 5.2  and Corollary 5.3 in \cite{H-const-twisted-mod}
(see Section 2).
\epfv

\begin{rema}
{\rm Since by Theorem \ref{first-existence}, $\widehat{M}^{[g]}_{B}$ is not $0$, 
Theorem \ref{2.5th-existence} in particular shows that there are nonzero 
lower-bounded generalized $g$-twisted $V$-modules when $V$ is a M\"{o}bius vertex superalgebra.}
\end{rema}

We also need to discuss the existence of 
lowest weights and lowest weight subspaces of lower-bounded 
generalized $g$-twisted $V$-modules.

\begin{defn}
{\rm A complex number $n_{0}$ is called a {\it lowest weight} of a lower-bounded 
generalized $g$-twisted $V$-module $W$ if $W_{[n_{0}]}\ne 0$ and 
$W_{[n]}=0$ for $n\in \C$ satisfying
$\Re(n)<\Re(n_{0})$. 
When lowest weights exist, 
we shall call the subspace of $W$ spanned by all homogeneous elements whose weights 
are lowest weights the {\it lowest weight space of $W$}. In general,
we call the infimum of the real parts of the weights of elements of $W$ 
the {\it weight infimum of $W$}. }
\end{defn}

\begin{rema}
{\rm Note that in general, lowest weights might not exist.
Even if they exist, in general they might not be unique;
they can be differed by imaginary numbers. In particular, the imaginary parts of the weights of the
elements of the lowest weight space can be different. }
\end{rema}

We need a result on the existence of the lowest weights of irreducible 
lower-bounded generalized $g$-twisted $V$-module. Recall that $P(V)$ is the set 
of $\alpha\in [0, 1)+\R$ such that $e^{2\pi i\alpha}$ is an eigenvalue of $g$ (on $V$). 

\begin{prop}\label{lowest-weight}
Let $V$ be a M\"{o}bius vertex superalgebra and $g$ an automorphism of $V$. Assume that the set
of real parts of the numbers in $P(V)$ has no accumulation point in $\R$. 
Then a finitely-generated lower-bounded generalized $g$-twisted $V$-module has a lowest weight.
In particular, an irreducible lower-bounded generalized $g$-twisted $V$-module has a lowest weight. 
\end{prop}
\pf
We need only prove that a lower-bounded generalized $g$-twisted $V$-module $W$ generated
by one homogeneous element $w$ has a lowest weight. 
By Definition \ref{generated-mobius},
$W$ is spanned by the coefficients of formal series of the form 
$Y^{g}_{W}(v, x)
L_{W}(-1)^{p}L_{W}(0)^{q}L_{W}(1)^{r}w$ for $v\in V$, $p, q, r\in \N$. When 
$v\in V^{[\alpha]}$, 
$$Y^{g}_{W}(v, x)=\sum_{s\in \Z}\sum_{k\in \N}
(Y^{g}_{W})_{\alpha+s, k}(v)x^{-\alpha-s-1}(\log x)^{k}.$$
So $W$ is spanned by elements of the form
\begin{equation}\label{element-form-4}
(Y^{g}_{W})_{\alpha+s, k}(v)L_{W}(-1)^{p}L_{W}(0)^{q}
L_{W}(1)^{r}w
\end{equation}
for $v\in V^{[\alpha]}$, $\alpha\in P(V)$,
$s\in \Z$, $k, p, q, r\in \N$. When $v$ is also homogeneous,
the weight of such an element is $\wt v-\alpha-s-1+p-r+\wt w$.
Since
$\wt v-\alpha-s-1+p-r\in -(P(V)+\Z)$, the real part of $\wt v-\alpha-s-1+p-r+\wt w$
is in $-(\Re(P(V))+\Z)+\Re(\wt w)$. 

Let $S$ be the set of 
the real parts of the weights of the nonzero homogeneous subspaces of $W$. 
Since $W$ is lower bounded, there exists a weight infimum $h\in \R$ of $W$, which by definition is
the infimum of $S$. Since $\wt v-\Re(\alpha)-s-1+p-r+\Re(\wt w)\ge h$, 
$\wt v-\Re(\alpha)-s-1+p-r\ge h-\Re(\wt w)$. In any neighborhood of $h$,
there is at least one element of $W$ such that the real part of its weight is
in the neighborhood. 
 Let $n\in \Z$ satisfying $n-1<h-\Re(\wt w)\le n$. Then for any nonzero 
element of the form 
(\ref{element-form-4}) of weight 
$\wt v-\alpha-s-1+p-r+\wt w$ such that $\wt v-\Re(\alpha)-s-1+p-r\le n$,
we also have $n-1<\wt v-\Re(\alpha)-s-1+p-r$. So we obtain
$$n-1-\wt v+s+1-p+r<-\Re(\alpha)\le n-\wt v+s+1-p+r.$$ 
Since $-1<-\Re(\alpha)\le 0$,
$n-1<\wt v-s-1+p-r<n+1$. But $\wt v-s-1+p-r$ is an integer. So we obtain $\wt v-s-1+p-r=n$.
Thus the weight of (\ref{element-form-4}) in this case is $n-\alpha+\wt w$.
Since $h$ is the infimum of $S$, there must be a sequence of elements of $W$ of the form 
(\ref{element-form-4}) such that limit of the real parts of the weights of the elements in the sequence
is $h$. In particular, we can assume that the real parts of the weights of the elements in the sequence 
is less than or equal to $n$. Then these real parts must be of the form $n-\Re(\alpha_{m})
+\Re(\wt w)$ for $m\in \Z_{+}$ and $\alpha_{m}\in P(V)$. Since these real parts form a sequence 
whose limit is $h$, the limit of the sequence $\{-\Re(\alpha_{m})\}_{m\in \Z_{+}}$ exists and is equal to 
$h-n-\Re(\wt w)$. If there are infinitely many different $-\Re(\alpha_{m})$ in the sequence, the limit of 
the sequence is an accumulation point of the set of the real parts of $P(V)$. But
 by our assumption, this set has no accumulation point in $\R$. 
Contradiction. So when $m$ is sufficiently large, $-\Re(\alpha_{m})$ are all equal and must be equal to the limit
$h-n-\Re(\wt w)$ of the sequence $\{-\Re(\alpha_{m})\}_{m\in \Z_{+}}$. 
Thus $h=n-\Re(\alpha_{m})+\wt w$ for $m$ sufficiently large. This proves that 
$n-\alpha_{m}+\wt w$ are lowest weights of $W$ for $m$ sufficiently large.
\epfv

\begin{exam}\label{finite-order-lw}
{\it Let $V$ be a M\"{o}bius vertex superalgebra and $g$ an automorphism of $V$
such that there are only finitely many distinct eigenvalues of $g$. 
Then $P(V)$ is a finite set. In particular, the set of the real parts of the elements of 
$P(V)$ does not have an accumulation point. 
So by Proposition \ref{lowest-weight}, in this case, a finitely-generated lower-bounded generalized $g$-twisted $V$-module
has a lowest weight. 
When the semisimple part $e^{2\pi i\mathcal{S}_{g}}$ of $g$ is of finite order (including the case that 
$g$ is of finite order), there are only finitely many distinct eigenvalues of $g$. 
Thus by Proposition \ref{lowest-weight}, 
a finitely-generated lower-bounded generalized $g$-twisted $V$-module
has a lowest weight when the semisimple part $e^{2\pi i\mathcal{S}_{g}}$ 
of $g$ is of finite order (in particular, when 
$g$ is of finite order). We see that Proposition \ref{lowest-weight} is indeed a generalization of the 
well-known fact that a finitely-generated lower-bounded generalized $g$-twisted $V$-module
has a lowest weight when $g$ is of finite order. }
\end{exam}

Let $W$ be a lower-bounded generalized $g$-twisted $V$-module. 
Then by Proposition 3.3 in \cite{H-twisted-int}, $W'$ has a structure of a
lower-bounded generalized $g^{-1}$-twisted $V$-module called the contragredient 
module of $W$. We have a functor $'$ called the contragredient functor. 

\begin{lemma}\label{sets-lowest-weights}
Let $V$ be a M\"{o}bius vertex superalgebra and $g$ an automorphism of $V$. Assume that the set
of real parts of the numbers in $P(V)$ has no accumulation point in $\R$. 
Then the set of the lowest weights of the irreducible  lower-bounded generalized $g$-twisted $V$-modules
and the set of the lowest weights of the irreducible lower-bounded generalized 
$g^{-1}$-twisted $V$-modules are the same.
\end{lemma}
\pf
By Proposition \ref{lowest-weight},
an  irreducible  lower-bounded generalized $g$-twisted $V$-module $W$ has 
lowest weights. Given a lowest weight $h$ of $W$, the dual space of the subspace of $W$ of weight $h$
 is the  subspace of $W'$ of the same weight $h$. Let $W_{0}$ be the submodule of $W'$ generated by an element 
$w'$ of  this subspace of $W'$. By Remark \ref{second-existence-mobius} and 
Theorem \ref{second-existence},
there exists a maximal submodule $J$ of $W_{0}$ such that $J$ does not contain $w'$
and the quotient $W_{0}/J$ is irreducible. It is clear that a lowest weight of $W_{0}/J$ 
is $h$, the weight of $w'$. But $h$ is a lowest weight of $W$. Therefore the set of 
 the lowest weights of the irreducible  lower-bounded generalized $g$-twisted $V$-modules
is contained in the set of the lowest weights of the irreducible lower-bounded generalized 
$g^{-1}$-twisted $V$-modules. Similarly, we see that 
the set of  the lowest weights of the irreducible  lower-bounded generalized 
$g^{-1}$-twisted $V$-modules
is contained in the set of the lowest weights of the irreducible lower-bounded generalized 
$g$-twisted $V$-modules. Thus these two sets are the same. 
\epfv

We are  now ready to study the existence of irreducible grading-restricted 
generalized or ordinary $g$-twisted $V$-module.

\begin{thm}\label{third-existence} 
Let $V$ be a M\"{o}bius vertex superalgebra and $g$ an automorphism of $V$.  Assume that the set
of real parts of the numbers in $P(V)$ has no accumulation point in $\R$. 
Also assume that the set of the real parts of the lowest weights of the irreducible
lower-bounded generalized $g$-twisted $V$-modules has a maximum such that 
the lowest weight subspace of an irreducible lower-bounded generalized $g$-twisted $V$-module
with this maximum as the real part of its lowest weight is finite dimensional.
Then 
there exists an irreducible grading-restricted generalized $g$-twisted $V$-module.
Such an irreducible grading-restricted generalized
$g$-twisted  $V$-module is an irreducible ordinary $g$-twisted  $V$-module if $g$ acts 
on it semisimply. In particular, if $g$ is of finite order, there exists an irreducible 
ordinary $g$-twisted  $V$-module.
\end{thm}
\pf
By Proposition \ref{lowest-weight}, every irreducible
lower-bounded generalized $g$-twisted $V$-module has a lowest weight. 
Let $W$ be an irreducible 
lower-bounded generalized $g$-twisted $V$-module such that the real part of 
its  lowest weights is the maximum of the set of the real parts of the lowest weights of the irreducible
lower-bounded generalized $g$-twisted $V$-modules and such that the lowest weight 
subspace of $W$ is finite dimensional. Then $W'$ also has a lowest weight
equal to the lowest weight of $W$ and the lowest weight subspace of $W'$ is also finite 
dimensional. 
By Lemma \ref{sets-lowest-weights},
the real part of the lowest weights
of $W'$  is  equal to the maximum of the set of
the real parts of the lowest weights of 
the irreducible lower-bounded generalized $g^{-1}$-twisted $V$-modules. 
If $W'$ is not generated by its lowest 
weight subspace, then the quotient of $W'$ by the submodule of $W'$ generated by 
its lowest weight space must be a lower-bounded generalized $g^{-1}$-twisted $V$-module 
such that the real part of the lowest weight is larger than the 
maximum above. But by Theorem \ref{second-existence}, the real part of the lowest weights of
this quotient of $W'$ must be equal to the real part of the lowest weights of an irreducible 
lower-bounded generalized $g^{-1}$-twisted $V$-module.  Contradiction. 
Thus $W'$ is generated by its lowest 
weight subspace. Thus $W'$ is a quotient of  $\widehat{M}^{[g]}_{B}$
with $M$ being its lowest weight subspace and with $B$ being 
its lowest weight. Since $M$ is finite dimensional and $\widehat{M}^{[g]}_{B}$ 
is spanned by elements of the form (\ref{element-form}),  the homogeneous subspaces of $\widehat{M}^{[g]}_{B}$
are of at most  countable dimensions. Thus homogeneous subspaces of $W'$ must 
also be of at most countable dimensions. This is possible only when 
the homogeneous subspaces of $W$ are finite dimensional. Thus $W$ is 
an irreducible grading-restricted generalized $g$-twisted $V$-module.

By Proposition \ref{semisimply}, such a grading-restricted generalized
$g$-twisted  $V$-module is an ordinary $g$-twisted  $V$-module if $g$ acts 
on it semisimply. 
\epfv

\begin{rema}
{\rm We assume in Theorem \ref{third-existence} 
that the set
of real parts of the numbers in $P(V)$ has no accumulation point in $\R$
  to make sure that every irreducible 
lower-bounded generalized $g$-twisted $V$-module has a lowest weight.
We can certainly replace this assumption by the weaker assumption that every irreducible 
lower-bounded generalized $g$-twisted $V$-module has a lowest weight.}
\end{rema}

Dong, Li and Mason proved in \cite{DLM1} (Theorem 9.1) that if $g$ is of finite order and 
the twisted Zhu's algebra $A_{g}(V)$ is not $0$ and finite dimensional, then
there exists an irreducible ordinary $g$-twisted $V$-module.
In Theorem 9.1 in \cite{DLM2}, Dong, Li and Mason proved that
if $V$ is a $C_{2}$-cofinite 
simple vertex operator algebra and $g$ is of finite order, then
there exists an irreducible ordinary $g$-twisted $V$-module. 
In the following consequence of Theorem \ref{third-existence},
we are able to remove the conditions that $A_{g}$ is not $0$ in 
Theorem 9.1 in  \cite{DLM1}, the simplicity of $V$ and
$C_{2}$-cofiniteness in Theorem 9.1 in \cite{DLM2} and 
weaken the conditions that $V$ has a conformal vector and 
$g$ is of finite order in these results in \cite{DLM1} and \cite{DLM2}:

\begin{cor}\label{fourth-existence}
Let $V$ be a M\"{o}bius vertex superalgebra and $g$ an automorphism of $V$.  Assume that the set
of real parts of the numbers in $P(V)$ has no accumulation point in $\R$. 
If $Z_{g}(V)$
(or, equivalently, $A_{g}(V)$) is finite dimensional, then there exists an irreducible grading-restricted generalized
$g$-twisted  $V$-module. Such an irreducible grading-restricted generalized
$g$-twisted  $V$-module is an irreducible ordinary $g$-twisted  $V$-module if $g$ acts 
on it semisimply. In particular, if $g$ is of finite order, there exists an irreducible 
ordinary $g$-twisted  $V$-module.
\end{cor}
\pf
In this case, there are only finitely many inequivalent irreducible $Z_{g}(V)$-modules. 
Moreover, these $Z_{g}(V)$-modules are all finite dimensional. 
Then using the functor $H_{g}$ in \cite{HY}
from the category of  lower-bounded $Z_{g}(V)$-modules
to the category of lower-bounded generalized $g$-twisted $V$-modules and 
Theorem \ref{second-existence}, we see that each irreducible $Z_{g}(V)$-module
gives an irreducible lower-bounded
generalized $g$-twisted $V$-module. By Proposition \ref{lowest-weight},
irreducible lower-bounded
generalized $g$-twisted $V$-modules have nonzero lowest weight spaces.
Taking the lowest weight spaces of 
irreducible lower-bounded
generalized $g$-twisted $V$-modules, we recover the 
irreducible $Z_{g}(V)$-modules that we start with. Thus we obtain a bijection between 
the set of equivalence classes of irreducible lower-bounded
generalized $g$-twisted $V$-modules and the set of 
irreducible  $Z_{g}(V)$-modules. Since  there are only finitely many 
inequivalent irreducible $Z_{g}(V)$-modules, there are also only finitely many 
inequivalent irreducible lower-bounded
generalized $g$-twisted $V$-modules. In particular, there exists a maximum of the 
real parts of the lowest weights of the irreducible  lower-bounded
generalized $g$-twisted $V$-modules. Also the lowest weight subspaces of 
irreducible  lower-bounded
generalized $g$-twisted $V$-modules are irreducible $Z_{g}(V)$-modules which are 
finite dimensional. By Theorem \ref{third-existence},
there exists a grading-restricted generalized 
$g$-twisted  $V$-module.
\epfv

\begin{rema}
{\rm In Theorem \ref{fourth-existence}, we assume that $V$ is a 
M\"{o}bius vertex superalgebra and 
that the set
of real parts of the numbers in $P(V)$ has no accumulation point in $\R$. 
But a vertex operator superalgebra is certainly a 
M\"{o}bius vertex superalgebra. Also by Remark \ref{finite-order-lw}, 
the assumption that $g$ is of finite order implies that 
the set of real parts of the numbers in $P(V)$ has no accumulation point in $\R$. 
Thus Theorem \ref{fourth-existence} is indeed a generalization of 
Theorem 9.1 in \cite{DLM1} and Theorem 9.1 in \cite{DLM2}.
}
\end{rema}

\begin{rema} 
{\rm Our proofs of Theorem \ref{third-existence}
and Corollary \ref{fourth-existence} do not need the results 
on genus-one $1$-point correlation functions. This fact shows that 
the existence of irreducible grading-restricted generalized or ordinary
$g$-twisted  $V$-module is not a genus-one property.}
\end{rema}

\renewcommand{\theequation}{\thesection.\arabic{equation}}
\renewcommand{\thethm}{\thesection.\arabic{thm}}
\setcounter{equation}{0}
\setcounter{thm}{0}
\section{Twisted extensions of modules for the fixed-point subalgebra}

In this section, using the lower-bounded generalized twisted module $M_{B}^{[g]}$
constructed in Section 5 of \cite{H-const-twisted-mod} (see Section 2),
we prove another result on the existence of certain lower-bounded generalized twisted modules. 

Let $V$ be a grading-restricted vertex algebra and $g$ an automorphism of $V$. 
Let $V^{g}$ be the fixed-point subalgebra of $V$ under $g$. 
Let $W$ be a lower-bounded generalized $g$-twisted $V$-module.
Then it is clear that $W$ equipped with the vertex operator map 
$Y_{W}^{g}\mbar_{V^{g}\otimes W}: V^{g}\otimes W\to W[[x, x^{-1}]]$ 
is a lower-bounded generalized $V^{g}$-module. Recall from \cite{H-twist-vo}
that $W=\coprod_{\alpha\in P_{W}}W^{[\alpha]}$ where 
$P_{W}$ is the set consisting of complex numbers $\alpha$ such that 
$\Re(\alpha)\in [0, 1)$ and $e^{2\pi i\alpha}$ is an eigenvalue 
of $g$ on $W$. Then we also know that $W^{[\alpha]}$ for $\alpha\in P_{W}$ 
are  lower-bounded generalized $V^{g}$-modules. 

One obvious question is: For a lower-bounded generalized $V^{g}$-module $W_{0}$, 
does  there exist  a lower-bounded generalized $g$-twisted $V$-module
$W$ such that  $W_{0}$ is a submodule of $W$ when $W$ is viewed as a 
lower-bounded generalized $V^{g}$-module? 
Unde the strong conditions on $V$ and the fixed point subalgebra $V^{G}$ under 
a finite group of automorphisms of $V$ that both $V$ and $V^{G}$ are $C_{2}$-cofinite and reductive in the sense
that every $\N$-gradable modules is a direct sum of irreducible module and 
all irreducible $g$-twisted $V$-modules (except for $V$) have positive weights, Dong, Ren and Xu proved
in \cite{DRX} that every irreducible $V^{G}$-module can be viewed as a submodule of an irreducible 
$g$-twisted $V$-module for some $g\in G$. 
Note that  a lower-bounded generalized $g$-twisted $V$-module in our definition always 
has an action of $g$ which can be written as 
$g=e^{2\pi i\mathcal{L}_{g}}$ where $\mathcal{L}_{g}=\mathcal{S}_{g}+\mathcal{N}_{g}$
and $\mathcal{S}_{g}$
and $\mathcal{N}_{g}$ are  semisimple and nilpotent operators on $W$. 
Thus one necessary condition for $W_{0}$ to be a submodule of $W$ is that we should also 
have  actions of $g$, $\mathcal{L}_{g}$, $\mathcal{S}_{g}$
and $\mathcal{N}_{g}$ on $W_{0}$.
In this section we prove that for such a  lower-bounded generalized $V^{g}$-module
$W_{0}$, the answer to the question above is positive. Note that the only condition in our result below 
is the necessary condition above.

\begin{thm}
Let $W_{0}$ be a lower-bounded generalized $V^{g}$-module (in particular, 
$W_{0}$ has a lower-bounded grading by $\C$ (graded by weights) and 
a grading by $\Z_{2}$ (graded by fermion numbers)).
Assume that 
$g$ acts on $W_{0}$ and there are semisimple and nilpotent operators $\mathcal{S}_{g}$
and $\mathcal{N}_{g}$, respectively, on $W_{0}$ such that 
$g=e^{2\pi i\mathcal{L}_{g}}$ where $\mathcal{L}_{g}=\mathcal{S}_{g}+\mathcal{N}_{g}$. 
Then $W_{0}$ can be extended to a lower-bounded generalized $g$-twisted $V$-module,
that is,  there exists a lower-bounded generalized $g$-twisted $V$-module $W$ 
and an injective module map $f: W_{0}\to W$ of $V^{g}$-modules. 
\end{thm}
\pf
By assumption, $W_{0}$ is a direct sum of generalized eigenspaces of $g$.
So we can assume that $W_{0}$ is a generalized  eigenspace of $g$ with eigenvalue 
$e^{2\pi i\alpha}$. 
Let $B$ be the largest real number such that for nonzero 
homogeneous $w\in W_{0}$, $\Re(\wt w)\ge B$.
Let $M$ be a space  of homogeneous 
generators of $W_{0}$ invariant under the action of $L_{W_{0}}(0)$,
$L_{W_{0}}(0)_{S}$, $L_{W_{0}}(0)_{N}$, $g$, $\mathcal{S}_{g}$, $\mathcal{N}_{g}$.  
Let $\{w^{a}\}_{a\in A}$ 
be a basis of $M$ such that $L_{W_{0}}(0)_{N}w^{a}=w^{L_{W_{0}}(0)_{N}(a)}$, where 
either $L_{W_{0}}(0)_{N}(a)\in A$ or $w^{L_{W_{0}}(0)_{N}(a)}=0$.
By Proposition \ref{span-1-v-o},  $W_{0}$ is spanned by coefficients of the formal series 
of the form $Y_{W_{0}}(v, x)L_{W_{0}}(-1)^{k}w^{a}$ for $v\in V^{g}$, $k\in \N$
 and $a\in A$. 

From Section 5 of \cite{H-const-twisted-mod} (see Section 2), we have a  lower-bounded 
generalized $g$-twisted $V$-module $M_{B}^{[g]}$. In particular, $M_{B}^{[g]}$
is a lower-bounded generalized $V^{g}$-module. Let $W_{1}$ be the 
$V^{g}$-submodule of $M_{B}^{[g]}$ generated by $(\psi^{a}_{M_{B}^{[g]}})_{-1, 0}\one$
for $a\in A$. Then by Proposition \ref{span-1-v-o} again,
$W_{1}$ is spanned by the coefficients of the formal series of 
the form 
$$Y_{W_{1}}(v, x)L_{W_{1}}(-1)^{k}(\psi^{a}_{M_{B}^{[g]}})_{-1, 0}\one
=Y_{M_{B}^{[g]}}(v, x)L_{M_{B}^{[g]}}(-1)^{k}(\psi^{a}_{M_{B}^{[g]}})_{-1, 0}\one$$ 
for $v\in V^{g}$, $k\in \N$ and $a\in A$. 

We define a linear map $f_{0}: W_{1}\to W_{0}$ using
$$
f_{0}\left(Y_{W_{1}}(v, x)L_{W_{1}}(-1)^{k}(\psi^{a}_{M_{B}^{[g]}})_{-1, 0}\one\right)
=Y_{W_{0}}(v, x)L_{W_{0}}(-1)^{k}w^{a}
$$
for $v\in V^{g}$, $k\in \N$ and $a\in A$. 
We first have to show that 
$f_{0}$ is well defined since there are relations among elements of the form 
$Y_{W_{1}}(v, x)L_{W_{1}}(-1)^{k}
(\psi^{a}_{M_{B}^{[g]}})_{-1, 0}\one$. But from Theorem \ref{spanning-set-2},
the only relations among elements of the form (\ref{element-form-2}) 
are  all induced from the relations among $v$'s. Hence these relations must also hold for any 
$V^{g}$-module, in particular for elements of the form 
$Y_{W_{0}}(v, x)L_{W_{0}}(-1)^{k}w^{a}$.
Thus $f_{0}$ is well defined.

Clearly $f_{0}$ preserves gradings. By definition, 
\begin{align*}
f_{0}((\psi^{a}_{\widehat{M}^{[g]}_{B}})_{-1, 0}\one)
&=w^{a},\\
f_{0}(L_{W_{1}}(-1)^{k}(\psi^{a}_{\widehat{M}^{[g]}_{B}})_{-1, 0}\one)
&=L_{W_{0}}(-1)^{k}w^{a},\\
f_{0}(Y_{W_{1}}(v, x)w)
&=Y_{W_{0}}(v, x)f_{0}(w)
\end{align*}
for $a\in A$, $k\in \N$ and $w\in W_{1}$. 
Using these formulas, 
we see that $f_{0}$ is a module map from $W_{1}$ to $W_{0}$. 
Clearly, $f_{0}$ is surjective. Thus there exists a lower-bounded generalized $V^{g}$-submodule
$J(W_{1})$ of $W_{1}$ such that $W_{1}/J(W_{1})$ is equivalent to $W_{0}$. 

Now let $J$ be the lower-bounded generalized $g$-twisted $V$-submodule 
of $M_{B}^{[g]}$ generated by $J(W_{1})$. Then $W=M_{B}^{[g]}/J$ is a 
 lower-bounded generalized $g$-twisted $V$-module containing a lower-bounded generalized 
$V^{g}$-submodule
$W_{2}=\{w+J\;|\; 
w\in W_{1}\}$. We now prove that this lower-bounded generalized 
$V^{g}$-module $W_{2}$ is isomorphic to $W_{1}/J(W_{1})$.
We define a linear map $f_{1}: W_{2}\to W_{1}/J(W_{1})$ by $f_{1}(w+J)=
w+J(W_{1})$ for $w\in W_{1}$. But we first need to prove that $f_{1}$ is well defined. 
This is equivalent to prove $W_{1}\cap J\subset J(W_{1})$. 

Let $w\in W_{1}\cap J$. Since $J$ is the lower-bounded generalized 
$g$-twisted $V$-submodule 
of $M_{B}^{[g]}$ generated by $J(W_{1})$ and $L_{M_{B}^{[g]}}(-1)$ 
acts on $J$, by Proposition \ref{span-1-v-o},
we can take $w$ to be a coefficient of 
$Y_{M_{B}^{[g]}}^{g}(u, x)w_{1}$ for $u\in V$ and 
$w_{1}\in J(w_{1})\subset W_{1}$. But $w$ is also in $W_{1}$.
Note that $W_{1}$ is a $V^{g}$-submodule of $M_{B}^{[g]}$ generated by 
$(\psi^{a}_{M_{B}^{[g]}})_{-1, 0}\one$
for $a\in A$ and, by the assumption that $W_{0}$ 
is a generalized  eigenspace of $g$ with eigenvalue 
$e^{2\pi i\alpha}$, $(\psi^{a}_{M_{B}^{[g]}})_{-1, 0}\one$ for $a\in A$ are 
generalized eigenvalues of $g$ with eigenvalue $e^{2\pi i\alpha}$. So 
$e^{2\pi i\mathcal{S}_{g}}w=e^{2\pi i\alpha}w$. 
If $u\in V^{[\beta]}$ for $\beta\in P(V)$, by (3.10) in \cite{H-twist-vo},
coefficients of 
$Y_{M_{B}^{[g]}}^{g}(u, x)w_{1}$ are generalized eigenvectors of $g$ 
with eigenvalue $e^{2\pi i(\alpha+\beta)}$. So if we take $w$ to be a coefficient of 
$Y_{M_{B}^{[g]}}^{g}(u, x)w_{1}$, $u$ must be in $V^{[0]}$.

If $u\in V^{[0]}$ is not in $V^{g}$, then $\mathcal{N}_{g}u\ne 0$.
On the other hand, since $w_{1}\in J(W_{1})\subset W_{1}$, we can take 
$w_{1}$ to be a linear combination  of coefficients of 
$Y_{W_{1}}(v, x)L_{W_{1}}(-1)^{k}
(\psi^{a}_{M_{B}^{[g]}})_{-1, 0}\one$ for $v\in V^{g}$,  
with $k\in \N$ and $a \in A$ to be fixed. 
We know that for $a\in A$, there exists $K_{a}\in \N$ such that 
$\mathcal{N}_{g}^{K_{a}+1}(\psi^{a}_{M_{B}^{[g]}})_{-1, 0}\one=0$
but $\mathcal{N}_{g}^{K_{a}}(\psi^{a}_{M_{B}^{[g]}})_{-1, 0}\one\ne 0$. 
Using  (3.12) in \cite{H-twist-vo}, $\mathcal{N}_{g}v=0$ for $v\in V^{g}$ and 
the fact that $\mathcal{N}_{g}$ commutes with $L_{W_{1}}(-1)$ which in turn 
follows from  (3.12) in \cite{H-twist-vo} and the $L(-1)$-derivative property for $Y_{W_{1}}$,  
we obtain
\begin{align}\label{twisted-extension-1}
\mathcal{N}_{g}&Y_{W_{1}}(v, x)L_{W_{1}}(-1)^{k}
(\psi^{a}_{M_{B}^{[g]}})_{-1, 0}\one\nn
&=Y_{W_{1}}(\mathcal{N}_{g}v, x)L_{W_{1}}(-1)^{k}
(\psi^{a}_{M_{B}^{[g]}})_{-1, 0}\one+
Y_{W_{1}}(v, x)\mathcal{N}_{g}L_{W_{1}}(-1)^{k}
(\psi^{a}_{M_{B}^{[g]}})_{-1, 0}\one\nn
&=Y_{W_{1}}(v, x)L_{W_{1}}(-1)^{k}
\mathcal{N}_{g}(\psi^{a}_{M_{B}^{[g]}})_{-1, 0}\one.
\end{align}
Then by (\ref{twisted-extension-1}), 
$$\mathcal{N}_{g}^{K_{a}+1}Y_{W_{1}}(v, x)L_{W_{1}}(-1)^{k}
(\psi^{a}_{M_{B}^{[g]}})_{-1, 0}\one=0.$$
But there exists $v$ such that 
$$Y_{W_{1}}(v, x)L_{W_{1}}(-1)^{k}\mathcal{N}_{g}^{K}
(\psi^{a}_{M_{B}^{[g]}})_{-1, 0}\one\ne 0.$$
Thus there exist an element $w_{1}$ of $J(W_{1})$ of the form of 
a linear combination  of coefficients of 
$Y_{W_{1}}(v, x)L_{W_{1}}(-1)^{k}
(\psi^{a}_{M_{B}^{[g]}})_{-1, 0}\one$ for $v\in V^{g}$, with fixed $k\in \N$ and $a\in A$
such that $\mathcal{N}_{g}^{K_{a}+1}w_{1}=0$ but 
$\mathcal{N}_{g}^{K_{a}}w_{1}\ne 0$. 

Now for $u\in V^{[0]}$ and $w_{1}$ of the form above,
by (3.12) in \cite{H-twist-vo}, 
\begin{align*}
\mathcal{N}_{g}Y_{M_{B}^{[g]}}^{g}(u, x)w_{1}
=Y_{M_{B}^{[g]}}^{g}(\mathcal{N}_{g}u, x)w_{1}+
Y_{M_{B}^{[g]}}^{g}(u, x)\mathcal{N}_{g}w_{1}.
\end{align*}
Then 
\begin{align}\label{twisted-extension-2}
\mathcal{N}_{g}^{K_{a}+1}Y_{M_{B}^{[g]}}^{g}(u, x)w_{1}
=\sum_{i=0}^{K_{a}+1}
\binom{K_{a}+1}{i}
Y_{M_{B}^{[g]}}^{g}(\mathcal{N}_{g}^{K_{a}+1-i}u, x)
\mathcal{N}_{g}^{i}w_{1}.
\end{align}
If there exists $u\in V^{[0]}$ such that 
$\mathcal{N}_{g}u\ne 0$, then we can always find $K_{u}\in \N$ such that 
$\mathcal{N}_{g}(\mathcal{N}_{g}^{K_{u}}u)=0$ but $\mathcal{N}_{g}^{K_{u}}u\ne 0$.
Taking $u$ to be $\mathcal{N}_{g}^{K_{u}}u$ if $K_{u}\ne 0$,
we see that we can always find $u\in V^{[0]}$ such that 
$\mathcal{N}_{g}^{2}u=0$ but $\mathcal{N}_{g}u\ne 0$. 
Then the right-hand side of (\ref{twisted-extension-2}) become 
$Y_{M_{B}^{[g]}}^{g}(\mathcal{N}_{g}u, x)
\mathcal{N}_{g}^{K_{a}}w_{1}$. From the construction 
of  $M_{B}^{[g]}$,   $Y_{M_{B}^{[g]}}^{g}(\mathcal{N}_{g}u, x)
\mathcal{N}_{g}^{K_{a}}w_{1}\ne 0$ if $\mathcal{N}_{g}u$ and 
$\mathcal{N}_{g}^{K_{a}}w_{1}$ are not $0$ because $M_{B}^{[g]}$ is universal. 
Since $\mathcal{N}_{g}u$ and 
$\mathcal{N}_{g}^{K_{a}}w_{1}$ are not $0$, $Y_{M_{B}^{[g]}}^{g}(\mathcal{N}_{g}u, x)
\mathcal{N}_{g}^{K_{a}}w_{1}$ cannot be $0$ in  $M_{B}^{[g]}$. 
So the right-hand side of  (\ref{twisted-extension-1})
is not $0$. In particular, there exist nonzero coefficients of the right-hand side of 
 (\ref{twisted-extension-1}). 
But when a coefficient of 
$Y_{M_{B}^{[g]}}^{g}(u, x)w_{1}$ is in  $W_{1}$, 
from our discussion above, we must have 
$$\mathcal{N}_{g}^{K_{a}+1}Y_{M_{B}^{[g]}}^{g}(u, x)w_{1}=0.$$
Contradiction. So  $\mathcal{N}_{g}u= 0$, that is, 
$u\in V^{g}$. Thus coefficients of $Y_{M_{B}^{[g]}}^{g}(u, x)w_{1}$
are in the generalized 
$V^{g}$-module generated by elements of $J(W_{1})$. But this generalized 
$V^{g}$-module is exactly $J(W_{1})$ itself. This shows that $w\in W_{1}\cap J$
must be in $J(W_{1})$, that is, $W_{1}\cap J\subset J(W_{1})$. This shows that 
$f_{1}$ is well defined.  

Clearly, $f_{1}$ is injective and surjective and preserves the gradings. So it is a 
grading-preserving linear isomorphism. By definition, $f_{1}(Y_{W_{1}}(v, x)w+J)
=Y_{W_{1}}(v, x)w+J(W_{1})$. Thus $f_{1}$ is an equivalence. 

Now we  can view 
$W_{1}/J(W_{1})$ as a lower-bounded generalized 
$V^{g}$-module of $W$.  
Let $f: W_{0}\to W_{1}/J(W_{1})$ be an equivalence 
of $V^{g}$-modules. Then $f$ can be viewed as an injective module map from 
$W_{0}$ to $W=M_{B}^{[g]}/J$.
\epfv

\noindent {\small \sc Department of Mathematics, Rutgers University,
110 Frelinghuysen Rd., Piscataway, NJ 08854-8019}

\noindent {\em E-mail address}: yzhuang@math.rutgers.edu

\end{document}